# POINT ESTIMATION WITH EXPONENTIALLY TILTED EMPIRICAL LIKELIHOOD[1]

By Susanne M. Schennach

*University of Chicago*

Parameters defined via general estimating equations (GEE) can be estimated by maximizing the empirical likelihood (EL). Newey and Smith [*Econometrica* **72** (2004) 219–255] have recently shown that this EL estimator exhibits desirable higher-order asymptotic properties, namely, that its $O(n^{-1})$ bias is small and that bias-corrected EL is higher-order efficient. Although EL possesses these properties when the model is correctly specified, this paper shows that, in the presence of model misspecification, EL may cease to be root $n$ convergent when the functions defining the moment conditions are unbounded (even when their expectations are bounded). In contrast, the related exponential tilting (ET) estimator avoids this problem. This paper shows that the ET and EL estimators can be naturally combined to yield an estimator called exponentially tilted empirical likelihood (ETEL) exhibiting the same $O(n^{-1})$ bias and the same $O(n^{-2})$ variance as EL, while maintaining root $n$ convergence under model misspecification.

**1. Introduction.** Statistical models defined via general estimating equations (GEE) of the form $E[g(x, \theta)] = 0$, where $g(x, \theta)$ is a vector-valued nonlinear function of a random vector $x$ and a parameter vector $\theta$, are very common in statistics. In such models, the parameter vector $\theta$ is traditionally estimated using two-step efficient generalized method of moments estimators (GMM) [21]. Over the last two decades, various one-step alternatives to two-step GMM have been suggested. Perhaps the best known estimators of this class are the empirical likelihood (EL), exponential tilting (ET) and GMM with continuous updating (CU) estimators, which have been previously studied in the econometrics [22, 26, 27, 35, 47] and statistics [37, 45, 48, 49, 50, 53] literatures. While all of these alternative estimators

Received July 2004; revised December 2005.
[1]Supported in part by NSF Grant SES-02-14068.

of $\theta$ share the first-order efficiency of efficient two-step GMM, their one-step nature provides them with desirable properties not enjoyed by GMM. In addition to bypassing the arbitrariness in the choice of first-step estimate (since any consistent estimate of $\theta$ can, in principle, be used as a first step and lead to slightly different second-step estimates in finite samples), these one-step estimators are also invariant under general parameter-dependent linear transformations of the vector of moment conditions [30, 50] and possess superior higher-order asymptotic properties [27, 28, 29, 47].

Considerable effort has been devoted to identifying which of these alternative estimators, EL, ET or CU, is preferable. Since all of these estimators are asymptotically equivalent up to $O_p(n^{-1/2})$ when the overidentifying restrictions are valid, differences must reside in their higher-order asymptotic properties or in their behavior under potential model misspecification. The CU estimator is generally regarded as less desirable than EL and ET because its objective function has often been observed to possess multiple modes [22, 30] and because it lacks the ability to generate likelihood ratio-based confidence regions whose shape adapts to the support of the data [4, 50]. Comparing ET and EL proves to be more difficult. On the one hand, based on a stochastic expansion argument, Newey and Smith [47] have established that EL should typically have a lower finite-sample bias than both ET and CU. Also, they have shown that bias-corrected EL is higher-order efficient than any other regular method of moments estimator. On the other hand, Imbens and co-workers [27, 30] have indicated that EL, unlike ET, exhibits a singularity in its influence function, suggesting that ET should be better behaved than EL in the presence of model misspecification. In addition, ET admits a computationally convenient treatment of misspecified models [32].

Although it can be argued that model misspecification can always be avoided through the use of specification tests, an alternative view is that most models are only approximations to the underlying phenomena and are therefore intrinsically misspecified. Accordingly, there exists a growing literature devoted to the study of so-called globally misspecified models (in which the misspecification does not vanish asymptotically). The classic theory of maximum likelihood estimators (MLE) when the distributional assumptions are misspecified can be found in [1, 25, 63, 64]. In this context, MLE consistently estimates the so-called *pseudo-true value* of the parameter of interest [56], which is defined as the parameter value associated with the distribution which is the closest to the true data generating process according to the so-called Kullback–Leibler information criterion (KLIC) discrepancy.

In recent years, the analysis of misspecified models has been actively extended to various extremum estimators [2, 13, 44, 51] and, in particular, to overidentified moment condition models [8, 18, 26, 32, 34, 41]. Overidentified models arise naturally in a number of applications. For instance, consider a regression model $y = x'\theta + \varepsilon$ where $\varepsilon$ is correlated with $x$ (so that least squares cannot be used) but uncorrelated with a vector of so-called instruments (denoted $z$). This leads to a vector of restrictions of the form $E[(y - x'\theta)z] = 0$, the dimension of which typically exceeds the dimension of $\theta$. Given the overidentified (i.e., overdetermined) nature of the restrictions, it is then possible that no value of $\theta$ simultaneously satisfies all the moment restrictions exactly in the population, resulting in a misspecified model [41]. A more extensive discussion of misspecified models as well as many references to empirical studies that perform inference with models which fail standard specification tests can be found in [18].

The motivation behind this interest for misspecified models stems from two observations. First, the imperfections of a model, although statistically detectable, may nevertheless be small in absolute terms and consequently have little impact on the results ([42], pages 1168–1169). Second, a misspecified but parsimoniously parametrized model may have better predictive power than a more realistic complex model which passes all specification tests ([9], page 596). At fixed sample size, as the number of parameters increases, their variances tend to increase as well, while the power of overidentification tests tends to decrease.

This paper is organized as follows. After briefly reviewing the properties of the EL, ET and CU estimators, we present a simple result that characterizes EL's poor behavior under misspecification in order to motivate the need for a new estimator. We then introduce an estimator called exponentially tilted empirical likelihood (ETEL) that naturally combines EL and ET, extending an approach previously considered in [10, 31, 40] for constructing likelihood-ratio confidence regions for the mean to the case of point estimation of parameters defined via general moment restrictions. The ETEL estimator is shown to be well behaved under model misspecification, like ET, while preserving the desirable higher-order asymptotic properties of EL established in [47]. Finally, simulations are used to illustrate the usefulness of this estimator. All proofs can be found in the Appendix.

## 2. Existing one-step alternatives to GMM.

### 2.1. *Generalities.* We first introduce our notation.

DEFINITION 1. Let $\theta$ denote the parameter vector of interest belonging to a compact subset $\Theta$ of $\mathbb{R}^{N_\theta}$. Let $x_i$ be sequence of random vectors taking values in $\mathcal{X} \subset \mathbb{R}^{N_x}$. Let $g(x_i, \theta)$ denote a vector of moment functions taking value in $\mathbb{R}^{N_g}$ and satisfying $E[g(x_i, \theta^*)] = 0$ at $\theta^* \in \Theta$. Let $n$ denote sample size and let all summations be over $1, \ldots, n$. Let $\|\cdot\|$ denote any convenient vector or matrix norm.

TABLE 1
*The EL, ET and CU estimators as particular cases of MED and GEL estimators (adapted from [47])*

| **Estimator** | $\gamma$ | $h(w)$ | $\rho(\xi)$ | $\tau(\xi)$ |
|---|---|---|---|---|
| EL | $-1$ | $-\ln nw$ | $\ln(1-\xi)$ | $(1-\xi)^{-1}$ |
| ET | $0$ | $nw \ln nw$ | $-\exp(\xi)$ | $-\exp(\xi)$ |
| CU | $1$ | $(nw)^2$ | $-(1+\xi)^2/2$ | $-(1+\xi)$ |
| ECR | $\gamma$ | $\frac{(nw)^{\gamma+1}-1}{\gamma(\gamma+1)}$ | $-\frac{1}{\gamma+1}(1+\gamma\xi)^{(\gamma+1)/\gamma}$ | $-(1+\gamma\xi)^{1/\gamma}$ |

The simplest way to summarize the properties of the EL, ET and CU estimators is to embed them in more general families of estimators. All three estimators admit two convenient representations. They can first be interpreted as minimum empirical discrepancy (MED) estimators [10, 11],

$$\hat{\theta} = \underset{\theta \in \Theta}{\arg\min} \left( n^{-1} \sum_i h(\hat{w}_i(\theta)) \right), \tag{1}$$

where $\hat{w}_i(\theta)$ is the solution to

$$\min_{\{w_i\}_{i=1}^n} n^{-1} \sum_i h(w_i) \tag{2}$$

subject to moment and normalization constraints,

$$\sum_i w_i g(x_i, \theta) = 0 \quad \text{and} \quad \sum_i w_i = 1. \tag{3}$$

(The term *empirical discrepancy* is used here to emphasize the fact that it is a discrepancy between measures supported on the sample rather than on a fixed discrete support.) Different choices of the discrepancy measure $h(\cdot)$ yield distinct estimators, as given in Table 1. Specific choices of $h(\cdot)$ have historically been given special names. The discrepancy used in EL, $h(w) = -\ln nw$, is known as the Kullback–Leibler information criterion (KLIC). Also, rewriting the minimization problem as an equivalent maximization problem, EL can be thought of as maximizing the "likelihood." In a similar fashion, ET, with $h(w) = nw \ln(nw)$, can be interpreted as maximizing a quantity known as *entropy*.

The minimum discrepancy formulation emphasizes that the estimator seeks to "reweight" the sample in order to satisfy the moment conditions exactly. The function $h(w_i)$ quantifies the amount of reweighting taking place and penalizes values that differ from $w_i = n^{-1}$. The point estimate $\hat{\theta}$ is the value that minimizes the discrepancy between $\hat{w}_i(\theta)$ and uniform weights. The weights $\hat{w}_i(\hat{\theta})$ are sometimes called *implied probabilities* because they can be used to construct more efficient empirical estimates of the data generating process [3, 26, 53].

EL, ET and CU can also be characterized as particular cases of the so-called generalized empirical likelihood (GEL) family of estimators [61],

$$
(4) \qquad \hat{\theta} = \underset{\theta \in \Theta}{\arg\min} \left( n^{-1} \sum_i \rho(\hat{\lambda}(\theta)' g(x_i, \theta)) \right),
$$

where the $N_g$-dimensional vector $\hat{\lambda}(\theta)$ is given by

$$
(5) \qquad \hat{\lambda}(\theta) = \underset{\lambda}{\arg\max} \left( n^{-1} \sum_i \rho(\lambda' g(x_i, \theta)) \right).
$$

The choice of the function $\rho(\cdot)$ defines the estimator used, as described in Table 1. The advantage of the GEL formulation is the computational convenience of solving an $(N_g + N_\theta)$-dimensional optimization problem rather than the $(n + N_\theta)$- dimensional problem defining an MED estimator.

As pointed out in [47], only specific choices of $h(\cdot)$ lead to an estimator admitting an equivalent GEL representation. A particularly rich class of such discrepancies is given by the Cressie–Read (CR) discrepancies [11],

$$
(6) \qquad h(w_i) = \frac{(nw_i)^{\gamma+1} - 1}{\gamma(\gamma+1)},
$$

where $\gamma$ is the parameter indexing the family. The corresponding $\rho(\cdot)$ is given in Table 1. The GEL representation of an MED estimator is called a dual problem because it amounts to reformulating the optimization in terms of the Lagrange multiplier $\lambda$ of the moment constraints. Newey and Smith [47] conjecture that Cressie–Read discrepancies may be the *only* discrepancies admitting a GEL representation. The weights attributed to the sample points in the original MED estimator can be recovered from

$$
(7) \qquad \hat{w}_i(\theta) = \frac{\tau(\hat{\lambda}(\theta)' g(x_i, \theta))}{\sum_j \tau(\hat{\lambda}(\theta)' g(x_j, \theta))},
$$

where $\tau(\xi) = d\rho(\xi)/d\xi$. We will refer to $\tau(\xi)$ as the tilting function because, as seen in equation (7), it indicates how the sample points are reweighted. The EL, ET and CU estimators are all members of this class (see Table 1) of empirical Cressie–Read (ECR) estimators. For a more detailed description of these families of estimators, we refer the reader to the excellent discussions found in [47, 50].

2.2. *Comparing the ECR estimators.* Let us first give the properties shared by all ECR estimators. For just-identified models ($N_g = N_\theta$), all of these estimators are trivially identical because the moment conditions can be satisfied exactly simply by choosing $\hat{\theta}$ appropriately without the need for tilting ($\hat{\lambda}(\hat{\theta}) = 0$). In over-identified models for which the over-identifying restrictions are valid, all ECR (and GEL) estimators possess the same asymptotic variance [47], which is equal to the asymptotic variance of the two-step efficient GMM estimator. All ECR estimators also enable the construction of confidence regions for the mean ($g(x_i, \theta) = x_i - \theta$) through convenient $\chi^2$-calibrated likelihood-ratio tests [4]. In light of the results in [53], Baggerly's results should extend to general $g(x_i, \theta)$.

The similarities end at the level of first order asymptotic properties in correctly specified models, however. As noted in [4], the behavior of the implied probabilities $\hat{w}_i(\hat{\theta})$ in finite samples differs markedly as a function of the sign of the parameter $\gamma$. For ECR with $\gamma \leq 0$, the implied probabilities $\hat{w}_i(\hat{\theta})$ are positive by construction, while for $\gamma > 0$, they can take on negative values. In a correctly specified model (where the implied probabilities converge to $n^{-1}$ for all ECR), negative weights become decreasingly likely as sample size grows and it is possible to entirely avoid negative weights via the use of a "shrinkage factor" correction (see [6]) that vanishes asymptotically and that has no effect on the limiting distribution. Nevertheless, under misspecification, the "shrinkage factor" correction does not vanish asymptotically since negative weights remain likely even asymptotically when $\gamma > 0$.

Positive implied probabilities are associated with likelihood-ratio confidence regions whose shape better adapts to the data [4, 50]. For instance, confidence regions for the mean then always lie within the convex hull of the support of the density of the corresponding random variable. Positive implied probabilities are also important in applications that require empirical estimates of the data generating process, as in the bootstrap, for instance, [7]. These observations indicate that EL (with $\gamma = -1$) and ET (with $\gamma = 0$) should be preferable to CU (with $\gamma = 1$). CU also suffers from a different problem, namely the potential presence of multiple local maxima in its objective function [22, 30].

Numerous authors have sought to further narrow down the choice of desirable ECR estimators. EL is often singled out among the ECR because it leads to likelihood ratio tests that are often, though not always, Bartlett correctable [10, 14, 39]. Newey and Smith [47] have recently shown that EL generally exhibits a smaller $O(n^{-1})$ bias than any other member of the ECR family [unless the centered third moments of the distribution of $g(x_i, \theta^*)$ happen to all vanish, in which case all ECR estimators have the same $O(n^{-1})$ bias]. They have also shown that bias-corrected EL is higher-order efficient, possessing an $O(n^{-2})$ variance that is no greater than that of any other bias-corrected regular method of moments estimator.

2.3. *Behavior under misspecification.* As mentioned in the Introduction, in the presence of misspecification, the object of interest is the pseudo-true value of the parameter vector. In the case of MED estimators, the pseudo-true value is defined as the value of $\theta$ which minimizes the population version of the empirical discrepancy used in the estimation procedure.

It is important to note that although two different estimators may consistently estimate the truth in a correctly specified model, they may converge in probability to different pseudo-true values in the presence of misspecification. These two pseudo-true values merely represent the minimizers of two different well-defined discrepancies. Even though it could be argued that pseudo-true values are generally "biased," the literature on estimation under model misspecification considers estimators of pseudo-true values as valid statistics for the purpose of inference (see [56], as an early reference). Following the recent literature using various ECR estimators under model misspecification, we will not argue whether any ECR has a "better" pseudo-true value than another in a given context. Instead, we will compare the convergence of various ECR estimators toward their respective pseudo-true values—a property that will be relevant regardless of the context of interest.

Imbens, Spady and Johnson [30] have informally argued that EL may be ill-behaved under model misspecification due to the fact that its influence function [20] is proportional to

$$\frac{1}{1 - \lambda' g(x_i, \theta^*)} \frac{\partial g'(x_i, \theta^*)}{\partial \theta} \lambda,$$

where the denominator $(1 - \lambda' g(x_i, \theta^*))$ can approach zero. We formalize this concern by showing that EL suffers from a dramatic degradation of its asymptotic properties under even the slightest amount of misspecification.

THEOREM 1. *Let $x_i$ be an i.i.d. sequence and assume $g(x, \theta)$ is twice continuously differentiable in $\theta$ for all $x$ and all $\theta \in \Theta$ and such that $\sup_{\theta \in \Theta} E[\|g(x_i, \theta)\|^2] < \infty$. If $\inf_{\theta \in \Theta} \|E[g(x_i, \theta)]\| \neq 0$ and $\sup_{x \in \mathcal{X}} u' g(x, \theta) = \infty$ for any $\theta \in \Theta$ and any unit vector $u$, then there exists no $\theta_{\mathrm{EL}}^* \in \Theta$ such that $\|\hat{\theta}_{\mathrm{EL}} - \theta_{\mathrm{EL}}^*\| = O_p(n^{-1/2})$.*

This theorem can be extended to the case where the moment function $g(x_i, \theta)$ diverges only along some directions $u$ but not others. In that case, the lack of root $n$ consistency is avoided only when $E[g(x_i, \theta^*)]$ happens to be orthogonal to the hyperplane along which $g(x_i, \theta)$ diverges.

While Theorem 1 does not prevent $\hat{\theta}_{\mathrm{EL}}$ from being a consistent estimator of the pseudo-true value $\theta_{\mathrm{EL}}^*$, it does preclude $\hat{\theta}_{\mathrm{EL}}$ from being root $n$ consistent, under the assumptions of the theorem. The proof of this result, which can be found in the Appendix, is somewhat involved, because standard asymptotics break down for EL under misspecification with unbounded $g(x, \theta)$. The following heuristic argument illustrates the nature of the problem: First note that the EL implied probabilities are given by

$$\hat{w}_i = n^{-1}(1 - \hat{\lambda}' g(x_i, \theta^*))^{-1} \tag{8}$$

and must be positive [49]. This implies that $\hat{\lambda} \xrightarrow{p} 0$, for otherwise, $\max_{i \leq n} \hat{\lambda}' g(x_i, \theta^*)$ would become unbounded as $n \to \infty$, causing some $\hat{w}_i$ to become negative. Now, the population version of the first order condition for $\hat{\lambda}$ is $E[g(x_i, \theta^*)/(1 - \lambda^{*'} g(x_i, \theta^*))] = 0$, where $\lambda^*$ and $\theta^*$ denote pseudo-true values. Yet, at the pseudo-true value $\lambda^* \equiv \text{plim}\, \hat{\lambda} = 0$, this expectation takes the value $E[g(x_i, \theta^*)]$, which is not zero, by the assumption of misspecification. Hence, the asymptotics of EL cannot be determined from a standard expansion of the first-order conditions around the pseudo-true values that satisfy the first-order conditions in the population. The limit as $n \to \infty$ and as $\hat{\lambda} \xrightarrow{p} 0$ cannot be freely exchanged, indicating that the moments entering the first-order conditions violate the standard dominance regularity conditions used to establish the asymptotics of $M$-estimators [46].

Theorem 1 indicates that, unless one is willing to solely use moment functions that take values in a compact set [so that $\sup_{x \in \mathcal{X}} u' g(x, \theta^*)$ is bounded for any $u$], the slightest amount of misspecification can cause the first-order asymptotic properties of EL to degrade catastrophically. It is important to realize that it is very common that the *function* $g(x, \theta)$ itself is unbounded even when $E[g(x, \theta)]$ is finite. For instance, if $g(x, \theta) = (x_1 - \theta, x_2 - 1)'$ and $x = (x_1, x_2)$ is drawn from a bivariate normal, $g(x, \theta)$ is unbounded even though $E[g(x, \theta)]$ exists.

Of course, when the main hypothesis of Theorem 1 ($\sup_{x \in \mathcal{X}} \|g(x, \theta)\| < \infty$) does not hold, root $n$ consistency becomes possible. For instance, the type of moment conditions advocated in the robustness literature (e.g., [20, 24]) involves bounded functions and root $n$ consistent estimation under misspecification, therefore possible using EL. Nevertheless, Theorem 1 does rule out moment conditions such as a simple average of random variables drawn from a distribution with unbounded support.

Theorem 1 is especially important given the growing literature on minimum empirical discrepancy estimators in misspecified models [8, 23, 26, 32, 34, 54, 61]. In the nonnested model selection literature using minimum discrepancies, it is often assumed that the competing models may be all misspecified and one is merely concerned with choosing the *least* misspecified model (e.g., [8, 32, 34]). Since the model that is eventually used for inference may then be misspecified, Theorem 1 is particularly relevant in this context and indicates that EL may not be well suited to these applications—unless the assumption of bounded $g(x_i, \theta)$ is made, which is precisely the assumption that the model selection literature using EL has so far relied upon [8, 23, 34].

EL's implied probability weights also exhibit questionable behavior under misspecification with unbounded $g(x_i, \theta)$. Since the EL implied probabilities $w_i = n^{-1}(1 - \hat{\lambda}' g(x_i, \theta))^{-1}$ must be positive [4], it is straightforward to see that $\hat{\lambda} \xrightarrow{p} 0$ when $g(x_i, \theta)$ is unbounded. Then note that the implied probabilities associated with all points $x_i$ such that $g(x_i, \hat{\theta}) \in C$ for a given compact set $C$ converge to $n^{-1}$ uniformly. Since this result holds for any compact set, this shows that, as sample size grows, all the adjustments to the implied probabilities become concentrated on the extreme observations. This would be desirable if the weights of these extreme observations were always decreased to ensure that the moment conditions are satisfied, but this is not the case. In fact, due to the convexity of EL's tilting function $\tau(\xi) = 1/(1-\xi)$, the reweighting of the sample in order to satisfy the misspecified moment conditions will be achieved by placing a large weight on a few extreme observations, while slightly reducing the weights (relative to $n^{-1}$) of the bulk of the observations. Note that this problem is exacerbated by the fact that the weights can become extremely large as the singularity in the tilting function is approached. This feature will be visible in our simulations below.

We conjecture that any ECR estimator with $\gamma < 0$ exhibits the same problems as EL under misspecification due to the presence of a ratio in the tilting function. Thus, if we focus solely on ECR which preclude negative implied probabilities ($\gamma \leq 0$), we are left with ET (corresponding to $\gamma = 0$) as the only candidate ECR whose behavior might not degrade dramatically under misspecification. This is precisely the choice made in [32] for the analysis of misspecified moment restriction models. The asymptotic variance of ET under misspecification is finite under reasonable assumptions, the most restrictive of which is slightly stronger than the requirement of the existence of the moment generating function $M_\theta(\lambda) = E[\exp(\lambda' g(x_i, \theta))]$ for $\theta$ and $\lambda$ in some bounded sets.

**3. Exponentially tilted empirical likelihood.** Higher-order asymptotic properties in correctly specified models point to EL, while good behavior under misspecification points toward ET. There appear to be significant benefits to be able to combine EL and ET into a single estimator exhibiting the advantages of both.

It has been suggested [10, 47, 50] that other GEL estimators that exhibit the same higher-order properties as EL can be devised by simply employing a tilting function $\tau(\xi)$ which admits the same Taylor expansion as the tilting function of EL in the vicinity of $\xi = 0$ up to sufficiently high order. The behavior of $\tau(\xi)$ farther away from $\xi = 0$ could then be independently set to match the behavior of ET under misspecification. This option is not particularly attractive because (i) the estimator completely loses its interpretation as a minimum empirical discrepancy estimator, (ii) the estimator can no longer be seen as either a maximum likelihood or a maximum entropy estimator, concepts that initially motivated the form of the EL and ET estimators, and (iii) there still exist an infinite number of ways to interpolate between EL and ET in order to construct $\tau(\xi)$, making the procedure highly nonunique. For these reasons we focus on a different approach.

3.1. *The estimator.* We propose to combine the EL and ET estimators in the following fashion.

DEFINITION 2 (*ETEL estimator*).

$$\hat{\theta} = \arg\min_{\theta}\left(n^{-1}\sum_i \tilde{h}(\hat{w}_i(\theta))\right), \tag{9}$$

where $\hat{w}_i(\theta)$ is the solution to

$$\min_{\{w_i\}_{i=1}^n} n^{-1}\sum_i h(w_i) \tag{10}$$

subject to

$$\sum_i w_i g(x_i,\theta) = 0 \quad \text{and} \quad \sum_i w_i = 1, \tag{11}$$

and where

$$\tilde{h}(w_i) = -\ln(nw_i), \tag{12}$$

$$h(w_i) = nw_i \ln(nw_i). \tag{13}$$

The discrepancies used in the above optimization problem correspond to using ET to find $\hat{w}_i(\theta)$ and using EL to find $\hat{\theta}$. Since $h(\cdot)$ belongs to the family of ECR discrepancies, this type of estimator still admits an $(N_g + N_\theta)$-dimensional dual optimization problem of the form

$$\hat{\theta} = \arg\min_{\theta} n^{-1}\sum_i \tilde{h}(\hat{w}_i(\theta)), \tag{14}$$

where $\hat{w}_i(\theta)$ is given by equation (7) with

$$\hat{\lambda}(\theta) = \arg\max_{\lambda}\left(n^{-1}\sum_i \rho(\lambda' g(x_i,\theta))\right) \tag{15}$$

and $\rho(\xi) = -\exp(\xi)$. This approach yields a unique estimator that combines the likelihood form of EL [equation (9)] while incorporating the concept of entropy characterizing ET [equation (10)]. For these reasons, we call this estimator exponentially tilted empirical likelihood (ETEL). Other authors [10, 31, 40] have considered this combination of EL and ET for the purpose of constructing likelihood-ratio confidence regions for the mean. It has also been shown that a nonparametric Bayesian procedure based on a prior on the space of distributions that favors distributions having a large entropy yields a posterior whose maximum would define the ETEL estimator [58]. This paper's contribution will be to identify the numerous desirable asymptotic properties of ETEL point estimates in the case of general moment functions $g(x_i, \theta)$ in the context of *overidentified* and possibly misspecified models.

The fact that the ETEL point estimate is the solution to two nested optimization problems (one of dimension $N_g$ and one of dimension $N_\theta$) instead of a single saddle-point problem does not complicate the implementation of the estimator. Indeed, ECR estimators are often implemented as two nested optimization problems despite their saddle point form, because it is easier to design robust numerical methods for locating either a maximum or a minimum that do not break down near inflection points of the objective function [43].

ETEL represents only one of the many possible combinations between two different discrepancies [one to find the $\hat{w}_i(\theta)$ and one to find $\hat{\theta}$]. However, using the EL discrepancy to find $\hat{\theta}$ stands out as a particularly attractive choice because the optimization problem defining $\hat{\theta}$ maintains the maximum likelihood form of EL, thus making it more likely that EL's higher-order properties will be preserved, an issue that will be investigated below. The use of the ET discrepancy to find the weights $\hat{w}_i(\theta)$ is also natural. Since the objective function for $\hat{\theta}$ contains $\ln(\hat{w}_i(\theta))$, it is imperative that the weights $\hat{w}_i(\theta)$ be positive by construction and not only asymptotically in correctly specified models. As noted earlier, if we focus on weights obtained from the ECR family, in order to maintain the low dimensional dual formulation, only ECR with $\gamma \leq 0$ provide positive weights by construction [4]. However, any ECR with $\gamma < 0$ contains a singularity in its influence function, leaving $\gamma = 0$, or ET, as the only sensible choice to find the weights in the presence of potential model misspecification.

From a conceptual point of view, one may wonder about the interpretation of the ETEL estimator, since its definition combines two different discrepancies. It is often pointed out that in the case of discrete distributions, EL provides maximum likelihood estimates of both $\theta^*$, the true value of the parameter vector of interest, *and* the weights. Since ET weights are used in ETEL, ETEL weights are not maximum likelihood estimates, but in itself this is not a great concern since the weights are nuisance parameters and inference focuses on $\theta^*$. Indeed, after solving for all the parameters in terms of $\theta$, both the ETEL and EL estimators of $\theta^*$ can be cast into the familiar form of a maximum likelihood estimator of $\theta^*$ (as opposed to both $\theta^*$ and the weights, as in EL),

$$\text{(16)} \qquad \hat{\theta} = \arg\max_\theta \left( n^{-1} \sum_i \ln(\hat{w}_i(\theta)) \right),$$

where $\hat{w}_i(\theta)$ is given by equations (7) and (5). Of course, such an estimator can only formally be identified as a maximum likelihood estimator in the special case of a discrete distribution having support consisting of a finite number of points. More generally, for continuous distributions we can nevertheless refer to $\hat{\theta}$ as a MED estimator of $\theta^*$ using the KLIC discrepancy (as for maximum likelihood estimators), an interpretation that will be relevant under model misspecification. The distinction between EL and ETEL lies in how the estimate of the distribution of the data generating process $\hat{w}_i(\theta)$ given $\theta$ is constructed. In a parametric likelihood, $\hat{w}_i(\theta)$ would be uniquely given by the distributional assumptions of the model. When moment conditions replace distributional assumptions, however, there exists no such unique choice of $\hat{w}_i(\theta)$, due to the nonparametric nature of the problem. Both EL and ETEL replace parametric distributional assumptions by a so-called *least favorable family* of distributions (see, e.g., [15]), that is, a parametric family of distributions (indexed by $\theta$) for which the estimation problem is as difficult as the original nonparametric problem. In other words, for each $\theta$ there exist an infinite number of distributions satisfying the moment conditions, and the specific discrete distribution defined by $\hat{w}_i(\theta)$ represents a worst-case scenario among them. As pointed out in [15], there exist numerous least favorable families; EL and ETEL merely employ different ones and, a priori, there is no reason to favor one over the other.

In the case of ETEL, the least favorable family chosen is the class of distributions obtained by maximizing entropy under the $\theta$-dependent moment constraints imposed by the model. Entropy maximization has a long history as a device to construct distributions which properly model lack of prior information under a set of known constraints (see, e.g., [12, 17, 36, 38, 60]). ETEL thus combines the well-established concept of entropy maximization to handle the nonparametric part of the estimation problem, while using likelihood maximization to deal with the parametric part of the problem. The idea of substituting nonparametric nonmaximum likelihood estimates [here, the $\hat{w}_i(\theta)$] into a likelihood-type objective function to avoid the pathological behavior of an approach based solely on maximum likelihood also parallels the work of Fan and co-workers [16].

One may have preferences regarding which estimator of the distribution is the more appealing, but the choice between EL and ETEL should ultimately be based on the comparison of the actual asymptotic properties of each estimator and their performance in simulation experiments, which is what the remainder of this article is devoted to.

3.2. *Properties.*

3.2.1. *First-order properties.* To simplify the notation, we make the dependence of all quantities on $\theta$ implicit and introduce the following definitions.

DEFINITION 3. Let $\hat{w}_i = \hat{w}_i(\theta)$, $\hat{\lambda} = \hat{\lambda}(\theta)$, $g_i = g(x_i, \theta)$, $\hat{g} = n^{-1} \times \sum_i g(x_i, \theta)$, $G_i = \partial g(x_i, \theta)/\partial \theta'$, $G = E[G_i]$, $\hat{G} = n^{-1} \sum_i G_i$, $\tilde{G} = \sum_i \hat{w}_i G_i$, $\hat{\Omega} = n^{-1} \sum_i g_i g_i'$, $\Omega = E[g_i g_i']$ and $\tilde{\Omega} = \sum_i \hat{w}_i g_i g_i'$. Quantities evaluated at $\theta = \theta^*$ are denoted by $*$.

Simple algebraic manipulations yield the following.

THEOREM 2. *The ETEL estimator $\hat{\theta}_{\text{ETEL}}$ maximizes the objective function*

$$\ln \hat{L}(\theta) = -\ln\left(n^{-1} \sum_i \exp(\hat{\lambda}'(g_i - \hat{g}))\right), \tag{17}$$

*where $\hat{\lambda}$ is such that*

$$n^{-1} \sum_i \exp(\hat{\lambda}' g_i) g_i = 0. \tag{18}$$

*The first-order conditions for $\hat{\theta}_{\text{ETEL}}$ can be written as*

$$n^{-1} \sum_i (1 - n\hat{w}_i) \frac{d(\hat{\lambda}' g_i)}{d\theta'} = 0, \tag{19}$$

*where the total derivative indicates that $\hat{\lambda}$ is allowed to vary with $\theta$.*

We then establish that ETEL is at least as good as any ECR estimator both in terms of its first-order asymptotic properties and in terms of its invariance properties.

ASSUMPTION 1 (*Regularity conditions*).

1. $x_i$ forms an i.i.d. sequence.
2. $\theta^* \in \text{int}(\Theta)$ is the unique solution to $E[g(x_i, \theta)] = 0$, where $\Theta$ is compact.
3. $g(x_i, \theta)$ is continuous (in $\theta$) at each $\theta \in \Theta$ with probability one.
4. $E[\sup_{\theta \in \Theta} \|g_i\|^{2+\delta}] < \infty$ for some $\delta > 0$ and $E[\sup_{\theta \in \mathcal{N}} \|G_i\|] < \infty$.
5. $\Omega^*$ is nonsingular and finite and $\text{rank}(G^*) = N_\theta$.
6. $g(x_i, \theta)$ is continuously differentiable (in $\theta$) in a neighborhood $\mathcal{N}$ of $\theta^*$.

These assumptions match those of Theorem 3.2 in [47] and include those of Theorem 3.4 in [50].

THEOREM 3 (First-order properties). *Under Assumption 1, the ETEL estimator (*i*) has the same limiting distribution as efficient two-step GMM,*

$$n^{1/2}(\hat{\theta} - \theta^*) \xrightarrow{d} N(0, \Sigma),$$

where $\Sigma = (G^{*\prime}(\Omega^*)^{-1}G^*)^{-1}$, (ii) ETEL enables the construction of $\chi^2$-calibrated likelihood-ratio confidence regions for $\theta$,

$$-2n\ln(\hat{L}(\theta)/\hat{L}(\hat{\theta})) \xrightarrow{d} \chi^2_{N_\theta},$$

and (iii) of $\chi^2$-calibrated test of the validity of overidentifying restrictions,

$$-2n\ln(\hat{L}(\hat{\theta})) \xrightarrow{d} \chi^2_{N_g-N_\theta}.$$

THEOREM 4 (Implied probabilities and invariance properties). *Whenever the ETEL estimator is defined, (i) it yields positive implied probabilities ($\hat{w}_i(\theta) \geq 0$), (ii) it is invariant under arbitrary one-to-one differentiable reparametrization $\theta = T(\beta)$ of the moment conditions [the estimate $\hat{\beta}$ obtained from the reparametrized moment conditions satisfies $\hat{\theta} = T(\hat{\beta})$] and (iii) it is invariant under general parameter-dependent nonsingular linear transformation $A(\theta)$ of the vector of moment conditions (using $E[A(\theta)g(x_i,\theta)] = 0$ or $E[g(x_i,\theta)] = 0$ as moment conditions gives the same $\hat{\theta}$).*

3.2.2. *Higher-order asymptotic properties.* Estimators having the same (first-order) asymptotic variance can be compared on the basis of their higher-order ($o_p(n^{-1/2})$) asymptotic properties [55]. While it has been established that likelihood-ratio confidence regions of the mean constructed using ETEL do not share EL's Bartlett correctability [10, 31], another type of analytic higher-order correction permits the same improvement in the order of the coverage accuracy [40]. Moreover, it has been observed in simulation studies [50, 62] that the Bartlett correction is often ineffective in practice because the "QQ" plots for the EL likelihood ratio test statistics are typically curved, making it unlikely that a linear correction such as Bartlett's would improve coverage accuracy. Finally, given that ETEL's objective function can be interpreted as a posterior for the parameter $\theta$ obtained via a nonparametric Bayesian procedure [58], it may be a more relevant and interesting topic of future research to verify whether a Bayesian Bartlett correction [5], which differs from the usual frequentist Bartlett correction, would be applicable to ETEL.

More importantly, we can show that the ETEL point estimate $\hat{\theta}_{\text{ETEL}}$ shares all of the other higher-order properties of EL established in [47]. Higher-order asymptotic properties of an estimator $\hat{\theta}$ are defined through a stochastic expansion (see, e.g., [52, 55]) of the form

(20) $$(\hat{\theta} - \theta^*) = n^{-1/2}\bar{\psi} + n^{-1}\bar{q} + n^{-3/2}\bar{r} + O_p(n^{-2}),$$

where $\bar{\psi}$, $\bar{q}$ and $\bar{r}$ are $O_p(1)$ and where $\bar{\psi}$ and $\bar{r}$ have zero mean. Within this framework, the $O(n^{-1})$ bias is defined as $E[\bar{q}]$ and represents the most important correction to standard first-order asymptotics based solely on the influence function $\bar{\psi}$. Another important correction to first-order asymptotics is the $O(n^{-2})$ variance, defined as

$$\text{(21)} \qquad \text{Var}[\bar{q}] + \text{Covar}[\bar{r}, \bar{\psi}] + \text{Covar}[\bar{\psi}, \bar{r}].$$

This expression can be informally obtained by computing the variance of equation (20). In general, it is not meaningful to compare the $O(n^{-2})$ variances of two estimators that possess different $O(n^{-1})$ biases and bias-corrected estimators should be used to compare efficiency.

We now proceed to compare the stochastic expansions of $\hat{\theta}_{\text{ETEL}}$ and $\hat{\theta}_{\text{EL}}$, using assumptions found in [47]. Our approach consists of establishing that the difference $\hat{\theta}_{\text{ETEL}} - \hat{\theta}_{\text{EL}}$ is such that the Newey and Smith results for $\hat{\theta}_{\text{EL}}$ carry over to $\hat{\theta}_{\text{ETEL}}$.

THEOREM 5 (Higher-order equivalence to EL). *Under Assumption 1 and if $E[\sup_{\theta \in \mathcal{N}} \|g_i\|^4] < \infty$, $E[\sup_{\theta \in \mathcal{N}} \|G_i\|^2] < \infty$ and for $\theta \in \mathcal{N}$, $G(x_i, \theta)$ is Lipschitz in $\theta$ with prefactor $b(x_i)$ such that $E[b(x_i)] < \infty$, then $\hat{\theta}_{\text{ETEL}} - \hat{\theta}_{\text{EL}} = O_p(n^{-3/2})$.*

A consequence of this result is that the ETEL estimator has the same $O(n^{-1})$ bias as the EL estimator obtained in [47], under their assumptions. (As shown in [59], this result in fact extends to all estimators constructed by substituting GEL weights into the EL objective function.)

ASSUMPTION 2. There exists a function $b(x_i)$ with $E[(b(x_i))^6] < \infty$ such that, in a neighborhood $\mathcal{N}$ of $\theta^*$, all partial derivatives of $g(x_i, \theta)$ with respect to $\theta$ up to order four exist, are bounded by $b(x_i)$ and are Lipschitz in $\theta$ with prefactor $b(x_i)$.

THEOREM 6 (Small bias property). *Under Assumptions 1 and 2, ETEL's $O(n^{-1})$ bias is*

$$n^{-1} H(-a + E[G_i H g_i])$$

*where $H = \Sigma G' \Omega^{-1}$ and $a$ is a vector whose elements are $a_j = \text{tr}(\Sigma E[\partial^2 g_j(x_i, \theta^*)/\partial\theta\,\partial\theta'])/2$, where $g_j(x_i, \theta^*)$ is the jth element of $g(x_i, \theta^*)$.*

Simple intuition for the small bias of EL is that the EL first-order condition resembles the first order condition of GMM ($\hat{g}'\hat{\Omega}^{-1}\hat{G} = 0$) except for the fact that the Hessian term $\hat{\Omega}$ and the Jacobian term $\hat{G}$ are replaced by efficient averages that are weighted by the EL implied probabilities [47]. This reweighting removes the $O(n^{-1})$ correlation between the different sample averages entering the first-order condition, thus reducing the bias. As shown in the Appendix, ETEL also efficiently weights the Hessian and the Jacobian terms, using only the ET weights. Since ET and EL implied probabilities are equivalent to a sufficiently high order, using the ET instead of the EL weights only contributes to a negligible $O_p(n^{-3/2})$ remainder.

The fact that ETEL and EL are equivalent up to $O_p(n^{-1})$ leads to two important simplifications in the comparison of their $O(n^{-2})$ variances. First, since their $O(n^{-1})$ biases are the same, the moments entering the expression for the bias correction of EL and ETEL are the same. If these moments were estimated in the same way for EL and for ETEL, then comparing the $O(n^{-2})$ variance of EL and ETEL with or without performing a bias correction would obviously give the same answer. This conclusion remains unchanged if the bias correction is applied using the EL estimate of $\theta$ for the EL bias correction and the ETEL estimate of $\theta$ for the ETEL bias correction since these estimates differ by $O_p(n^{-3/2})$, which would give rise to a difference of only $O_p(n^{-1}n^{-3/2})$ in the bias correction. Moreover, as pointed out in [47], whether the moments entering the bias correction are estimated by sample averages or averages weighted by implied probabilities has no effect on the higher-order variance of the resulting bias-corrected estimator. Hence, using EL weights for the EL bias correction and ETEL weights for the ETEL bias correction makes no difference either. In conclusion, we can meaningfully compare the $O(n^{-2})$ variances of EL and ETEL before performing a bias correction.

The second simplification made possible by the equivalence of the $O_p(n^{-1})$ terms of the EL and ETEL stochastic expansions is that the differences in their $O(n^{-2})$ variance must take the form

$$(22) \qquad \mathrm{Covar}[\bar{r}^{\mathrm{ETEL}} - \bar{r}^{\mathrm{EL}}, \bar{\psi}] + \mathrm{Covar}[\bar{\psi}, \bar{r}^{\mathrm{ETEL}} - \bar{r}^{\mathrm{EL}}],$$

as seen in (21). Hence, it is possible for ETEL and EL to differ by $O_p(n^{-3/2})$, while still sharing the same $O(n^{-2})$ variance, as long as that difference is uncorrelated with their (identical) influence function $\bar{\psi}$. In fact, this is precisely the case, as shown in the Appendix.

THEOREM 7 (Higher-order efficiency). *Under Assumptions* 1 *and* 2, *the $O(n^{-2})$ variances of ETEL and EL are equal.*

Maintaining the maximum likelihood form for the optimization problem defining $\hat{\theta}_{\mathrm{ETEL}}$ thus achieves the desired goal, namely, to maintain the higher-order asymptotic properties of EL found in [47]. It is the fact that (22) vanishes that enables ETEL to be higher-order efficient even though it differs sufficiently from EL to fail to be Bartlett correctable in the frequentist sense.

3.2.3. *Behavior under misspecification.* While in the previous section we have seen that ETEL inherits the higher-order properties of EL, we will now show that it also exhibits some of the desirable properties of ET that EL lacks under model misspecification.

Following the discussion of Section 3.1, ETEL's pseudo-true value $\theta^*$ minimizes the KLIC discrepancy between the true data generating process and an entropy maximizing least favorable family of distributions parametrized by $\theta$ (which replaces the distributional assumptions in parametric maximum likelihood).

We will now study the first-order asymptotic properties of ETEL under misspecification.

THEOREM 8. *For a given $\theta$, assume that $E[\exp(\lambda' g(x_i, \theta))]$ exists in a neighborhood of its minimum. If a subvector of $g(x_i, \theta)$ is statistically independent of the remaining elements of $g(x_i, \theta)$, then the empirical c.d.f. obtained from ETEL (or ET) implied probabilities at $\theta$ converges pointwise (at every point of continuity of the true c.d.f.) to a c.d.f. that maintains this independence, even under misspecification. EL achieves this only in the absence of misspecification.*

This indicates the possibility that using an empirical c.d.f. obtained from the implied probability weights of EL in the hope of improving accuracy could actually result in the introduction of a spurious dependence among variables. ETEL avoids this unappealing eventuality. This property could be helpful when the implied probabilities are employed to improve the efficiency of the bootstrap, as in [7], when the model happens to be misspecified.

A more important quality that ETEL shares with ET is the nonsingular behavior of its influence function. As noted by [30], an estimator's influence function $\psi(x_i)$ is proportional to its first-order conditions. By inspection of ETEL's first-order condition [equation (19)] it is clear ETEL's influence function will not contain any singularity, unlike EL's influence function. It will therefore not be surprising that ETEL avoids EL's undesirable behavior under misspecification, under regularity conditions similar to the ones made by [32] for ET, as shown more formally below.

Let $\lambda^*(\theta)$ denote the solution to $E[\exp(\lambda' g(x_i, \theta)) g(x_i, \theta)] = 0$, which is unique by the strict convexity of $E[\exp(\lambda' g(x_i, \theta))]$ in $\lambda$.

ASSUMPTION 3 (*Regularity conditions under misspecification*).

1. $x_i$ forms an i.i.d. sequence.
2. The function $\ln L(\theta) \equiv -\ln(E[\exp(\lambda^{*\prime}(\theta)(g(x_i, \theta) - E[g(x_i, \theta)]))])$ is maximized at a unique "pseudo-true" value $\theta^* \in \text{int}(\Theta)$, where $\Theta$ is compact.
3. $g(x_i, \theta)$ is continuous (in $\theta$) at each $\theta \in \Theta$ with probability one.

4. $E[\sup_{\theta \in \Theta} \sup_{\lambda \in \Lambda(\theta)} \exp(\lambda' g(x_i, \theta))] < \infty$ where $\Lambda(\theta)$ is a compact set such that $\lambda^*(\theta) \in \text{int}(\Lambda(\theta))$.
5. $S_{jl}(x_i, \theta) = \partial^2 g(x_i, \theta^*)/\partial \theta_j \partial \theta_l$ is continuous (in $\theta$) for $\theta \in \mathcal{N}$, a neighborhood of $\theta^*$.
6. There exists $b(x_i)$ satisfying $E[\sup_{\theta \in \mathcal{N}} \sup_{\lambda \in \Lambda(\theta)} \exp(k_1 \lambda' g(x_i, \theta)) \times (b(x_i))^{k_2}] < \infty$ for $k_1 = 1, 2$ and $k_2 = 0, 1, 2, 3, 4$ such that $\|G(x_i, \theta)\| \leq b(x_i)$ and $\|S_{jl}(x_i, \theta)\| \leq b(x_i)$ for $j, l = 1, \ldots, N_\theta$ for any $x_i \in \mathcal{X}$ and for all $\theta \in \mathcal{N}$.

The simplest way to describe the asymptotics of ETEL under misspecification is to introduce an equivalent just-identified GMM estimator involving an augmented parameter vector $\beta = (\tau, \kappa', \lambda', \theta')'$. The vector $\theta \in \mathbb{R}^{N_\theta}$ is the parameter vector of interest, while $(\tau, \kappa', \lambda')' \in \mathbb{R}^{1+2N_g}$ are auxiliary parameters to be estimated jointly with $\theta$. The dimension of this augmented parameter vector is higher than in the case of GEL estimators under misspecification ($1 + 2N_g + N_\theta$ instead of $N_g + N_\theta$). This is due to the fact that the first-order conditions for $\hat{\theta}$ in ETEL involve a few additional terms taking the form of a product of sample moments that are absent in GEL estimators. Each of these products of sample moments can be linearized by introducing the additional parameters $\kappa$ and $\tau$. Note that these additional parameters are merely a device used to simplify the construction of the covariance matrix of the estimator. The point estimate $\hat{\theta}$ can be obtained without introducing $\kappa$ and $\tau$, as seen in Theorem 2.

LEMMA 9. *The ETEL point estimate $\hat{\theta}$ is given by the appropriate subvector of the vector $\hat{\beta} = (\hat{\tau}, \hat{\kappa}', \hat{\lambda}', \hat{\theta}')'$, the solution to*

$$n^{-1} \sum_i \phi(x_i, \hat{\beta}) = 0,$$

*where, letting $\hat{\tau}_i = \exp(\hat{\lambda}' g_i)$,*

(23)
$$\phi(x_i, \hat{\beta}) = \begin{bmatrix} \hat{\tau}_i - \hat{\tau} \\ \dfrac{\partial}{\partial \hat{\kappa}}(\hat{\tau}_i g_i' \hat{\kappa} + \hat{\tau} g_i' \hat{\lambda} - \hat{\tau}_i) \\ \dfrac{\partial}{\partial \hat{\lambda}}(\hat{\tau}_i g_i' \hat{\kappa} + \hat{\tau} g_i' \hat{\lambda} - \hat{\tau}_i) \\ \dfrac{\partial}{\partial \hat{\theta}}(\hat{\tau}_i g_i' \hat{\kappa} + \hat{\tau} g_i' \hat{\lambda} - \hat{\tau}_i) \end{bmatrix}$$
$$= \begin{bmatrix} \hat{\tau}_i - \hat{\tau} \\ \hat{\tau}_i g_i \\ (\hat{\tau} - \hat{\tau}_i) g_i + \hat{\tau}_i g_i g_i' \hat{\kappa} \\ \hat{\tau}_i G_i' \hat{\kappa} + \hat{\tau}_i G_i' \hat{\lambda} g_i' \hat{\kappa} - \hat{\tau}_i G_i' \hat{\lambda} + \hat{\tau} G_i' \hat{\lambda} \end{bmatrix}.$$

TABLE 2
*The bias of the EL, ETEL and ET estimators for the Hall–Horowitz design*

|        | EL    | ETEL  | ET    |
|--------|-------|-------|-------|
| $K=4$  | 0.063 | 0.061 | 0.103 |
| $K=10$ | 0.129 | 0.103 | 0.232 |

Given the just-identified nature of the estimator defined in Lemma 9, its asymptotic distribution follows quite directly.

THEOREM 10 (Asymptotics under misspecification). *Let* $\Gamma = E[\partial \phi(x_i,\beta)/\partial \beta'|_{\beta=\beta^*}]$ *and* $\Phi = E[\phi(x_i,\beta^*)\phi'(x_i,\beta^*)]$. *Under Assumption* 3, *if* $\Gamma$ *is nonsingular, then* $n^{1/2}(\hat{\beta} - \beta^*) \xrightarrow{d} N(0, \Gamma^{-1}\Phi(\Gamma')^{-1})$.

**4. Simulations.** We first illustrate the fact that ETEL has the same $O(n^{-1})$ bias as EL. We use the simple experimental design suggested in [19] and subsequently used in [30, 33], slightly expanded to have $K$ moment conditions rather than two. The moment conditions are

$$g(x_i,\theta) = [r(x_i,\theta) \quad r(x_i,\theta)x_{i2} \quad r(x_i,\theta)(x_{i3}-1) \quad \cdots \quad r(x_i,\theta)(x_{iK}-1)]',$$

where $r(x_i,\theta) = \exp(-0.72 - (x_{i1}+x_{i2})\theta + 3x_{i2}) - 1$. These restrictions are satisfied at $\theta^* = 3$ when $(x_{i1},x_{i2})' \sim N(0,(0.16)I)$ and $x_{ik} \sim \chi_1^2$, for $k=3,\ldots,K$. Note that the third moments of all elements of $g(x_i,\theta)$ are nonzero and that $g(x_i,\theta)$ is nonlinear in $\theta$, so that the $O(n^{-1})$ bias does not trivially vanish. Figure 1 shows the c.d.f. of the EL, ET and ETEL estimators of $\theta$ obtained from 10,000 replicated samples of the above design (with $K=4$ and $K=10$), each containing 200 observations. [Samples for which at least one of the three estimators considered failed to converge were discarded. This happened 14 times for $K=4$ and 32 times for $K=10$. The most frequent reason for failure of convergence was that the origin was not contained within the convex hull of the values of $g(x_i,\theta)$ for any $\theta$, in which case none of the estimators is even defined. The number of nondiscarded samples is 10,000.] It is apparent that the ETEL and EL point estimates have very similar distributions, as expected from their equivalence up to the $O_p(n^{-1})$ term of their stochastic expansion. The distribution of the ET point estimates differs noticeably from that of EL and ETEL, and the main difference takes the form of a bias, which is reported in Table 2. The bias of ET increases more rapidly with the number of moment conditions than the biases of both EL and ETEL, as the higher-order asymptotics analyse given in [47] and in the present work would suggest.

Our next simulation compares the behavior of EL, ET and ETEL under misspecification. We consider a simple case where we wish to estimate the mean while imposing a known variance. In this example, the moment conditions are

$$g(x_i,\theta) = [\,x_i - \theta \quad (x_i - \theta)^2 - 1\,],$$

where $x_i$ is drawn either from a correctly specified Model C or a misspecified model M,

$$x_i \sim N(0,1) \qquad \text{(for Model C)},$$
$$x_i \sim N(0,(0.8)^2) \qquad \text{(for Model M)}.$$

Note that this experiment is specifically designed so that the pseudo-true value ($\theta^* = 0$) for the misspecified model is the same for EL, ET and ETEL, thus enabling a meaningful comparison of the variances of these estimators.

Figure 2 shows the c.d.f. of the EL, ET and ETEL estimators of $\theta$ for a sample size of 1000 and a sample size of 5000, evaluated with 10,000 and 2000 replications, respectively. The variability of the EL estimate is clearly larger than that of ET and ETEL, as confirmed by the calculated standard deviations given in Table 3. Interestingly, the distributions (and the standard deviations) of the ET and ETEL estimators are quite similar. While the ET and ETEL standard deviations shrink by the expected factor of $\sqrt{5}$ as the sample size is increased from 1000 to 5000, the standard deviation of EL barely changes, which is not surprising given the results of Theorem 1. Note that the difference between the distribution of EL and that of the two other estimators can be made arbitrarily large either by increasing the amount of misspecification or by increasing the sample size.

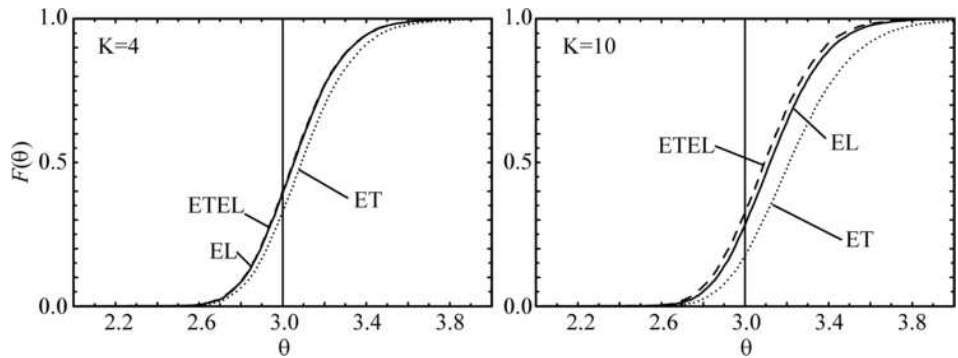

Fig. 1. *Cumulative distribution function of the EL, ET and ETEL estimators for the Hall–Horowitz design with* 4 (left) *and* 10 (right) *moment conditions. The sample size is* $n = 200$ *and* 10000 *replications were used to calculate this empirical c.d.f.*

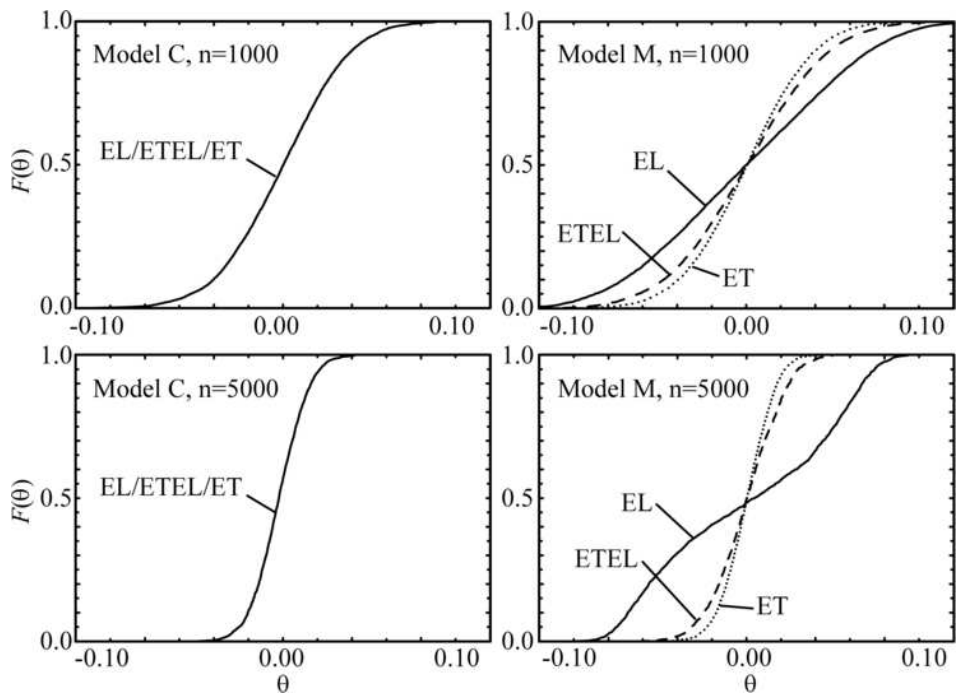

Fig. 2. *Cumulative distribution function of the EL, ET and ETEL estimators for Models C and M defined in the text. For the* top *portion of the figure, the sample size is* $n = 1000$ *and* 10000 *replications were used. For the* bottom *portion of the figure,* $n = 5000$ *and* 2000 *replications were used.*

We can also use simulations to illustrate the source of EL's poor behavior under misspecification. Figure 3 shows the implied probabilities for EL and ETEL in two simulated samples of size $n = 1000$ and $n = 5000$ drawn from the misspecified Model M. It is apparent that the EL implied probabilities attribute an excessive weight to the extreme observations. As the sample size grows, this trend worsens: the second graph exhibits an extremely large weight at $x_i \approx -3$ and $nw_i \approx 95$. In contrast, the ETEL implied probabilities distribute the weight more uniformly over the whole sample and, even more importantly, the weights do not become increasingly concentrated in the tails as the sample size grows.

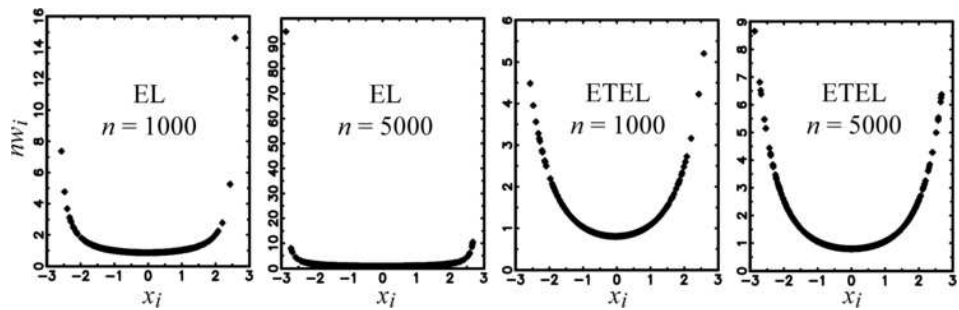

Fig. 3. *EL and ETEL implied probabilities in simulated samples drawn from the misspecified Model M as a function of sample size. Note the differences in the scale of the vertical axes.*

These examples, although simple and perhaps not realistic, illustrate how ETEL matches the low-bias property of the EL estimator and shares the reasonable behavior of ET under misspecification.

**5. Conclusion.** Our first important result is to show that although empirical likelihood (EL) is known to exhibit numerous desirable higher-order asymptotic properties in correctly specified models, its first-order asymptotic properties can degrade catastrophically in the presence of the slightest amount of misspecification, causing the loss of root $n$ consistency. Although the use of only bounded functions $g(x_i, \theta)$ in the moment conditions $E[g(x_i, \theta)] = 0$ avoids this problem, this is a rather strong constraint. In contrast, exponential tilting (ET) is known to be inferior to EL in terms of its higher-order properties, but remains well behaved in the presence of misspecification under relatively weak regularity conditions [32].

Our second main contribution is to show that EL and ET can be combined to yield an estimator that exhibits the advantages of both. This so-called exponentially tilted empirical likelihood (ETEL) has the same low $O(n^{-1})$ bias and the same $O(n^{-2})$ variance as EL in correctly specified models, and yet avoids EL's pitfalls in misspecified models.

## APPENDIX: PROOFS

The quantities given in Definitions 1 and 3 will be used throughout the Appendix. Let $C$ denote a generic constant which may take distinct values in different contexts. Let CSI stand for Cauchy–Schwarz inequality and let w.p.a. 1 stand for the phrase "with probability approaching one."

Table 3
*The standard deviations of the EL, ETEL and ET estimators for Models C and M defined in the text. The number of replications is* 10000 *for the* $n = 1000$ *sample and* 2000 *for the* $n = 5000$ *sample*

|           | $n = 1000$ |       |       | $n = 5000$ |       |       |
|-----------|------------|-------|-------|------------|-------|-------|
| Estimator | EL         | ETEL  | ET    | EL         | ETEL  | ET    |
| Model C   | 0.032      | 0.032 | 0.032 | 0.014      | 0.014 | 0.014 |
| Model M   | 0.054      | 0.038 | 0.031 | 0.052      | 0.019 | 0.014 |

PROOF OF THEOREM 1. The proof proceeds by constructing a triangular array of estimators $\hat{\theta}_{k,n}$ indexed by the sample size $n$ and by an auxiliary truncation parameter $k$. To define this array, let $\mathcal{G}_k$ be an increasing sequence of nested compact subsets of $\mathbb{R}^{N_g}$ such that $\bigcup_{k=1}^{\infty} \mathcal{G}_k = \mathbb{R}^{N_g}$. Then let $\mathcal{C}_k = \{x \in \mathcal{X} : g(x,\theta) \in \mathcal{G}_k \text{ for all } \theta \in \Theta\}$. Note that $\mathcal{C}_k$ is nonempty for $k$ sufficiently large.

Let $F_\infty(x)$ denote the distribution of $x$ and let $\hat{\theta}_{\infty,n}$ denote the EL estimator obtained from a sample of size $n$ and let $\theta^*_\infty$ denote EL's pseudo-true value, assuming it exists (for otherwise, $\hat{\theta}_{\infty,n}$ could not even be consistent).

Let $F_k(x)$ be a sequence of distributions indexed by $k \in \mathbb{N}$, each having support $\mathcal{C}_k$. We choose $F_k(x)$ so that, for all sufficiently large $k$, the moment conditions are uniformly misspecified ($\inf_{k \geq \bar{k}} \inf_{\theta \in \Theta} \|E_{F_k}[g(x,\theta)]\| > 0$ for some $\bar{k} \in \mathbb{N}$). Let $\hat{\theta}_{k,n}$ denote the EL estimator in a sample size of $n$ when the true data generating process is $F_k(x)$ and let $\theta^*_k \in \Theta$ denote the corresponding pseudo-true value. We then note that it is also always possible to choose a distribution $F_k(x)$ with support $\mathcal{C}_k$ such that $P[|u'(\hat{\theta}_{k,n} - \theta^*_k)| \geq \varepsilon] \leq P[|u'(\hat{\theta}_{\infty,n} - \theta^*_\infty)| \geq \varepsilon]$ for any $\varepsilon > 0$, any conformable unit vector $u$ and all $n$. For instance, one could first construct a distribution $\tilde{F}_k(x)$ equal to $F_\infty(x)$ conditional on the event $x \in \mathcal{C}_k$. Let $\theta^*_k$ denote the pseudo-true value associated with $\tilde{F}_k(x)$. Then set $F_k(x)$ to be a mixture of $\tilde{F}_k(x)$ and a degenerate distribution that would give $\theta^*_k$ as an EL estimate with certainty. In this fashion, $F_k(x)$ is a "truncated" version of $F_\infty(x)$ designed to make the estimation of $\theta^*_k$ by $\hat{\theta}_{k,n}$ easier than the estimation of $\theta^*_\infty$ by $\hat{\theta}_{\infty,n}$. Obviously, $\hat{\theta}_{k,n}$ is an infeasible estimator that uses out of sample information. It is introduced solely for the purpose of facilitating the proof. Note that $\theta^*_k \neq \theta^*_\infty$ in general, but the proof will never require that $\theta^*_k = \theta^*_\infty$.

For a distribution $F_k(x)$ having compact support, the EL estimator can be written as a just identified GMM estimator of an augmented parameter vector $\hat{\beta} = (\hat{\theta}'_{k,n}, \hat{\lambda}'_{k,n})'$ satisfying the first-order conditions

$$(24) \qquad n^{-1} \sum_i G'(x_i, \hat{\theta}_{k,n}) \hat{\lambda}_{k,n} / (1 - \hat{\lambda}' g(x_i, \hat{\theta}_{k,n})) = 0,$$

$$(25) \qquad n^{-1} \sum_i g(x_i, \hat{\theta}_{k,n}) / (1 - \hat{\lambda}' g(x_i, \hat{\theta}_{k,n})) = 0.$$

Note that these first-order conditions form a just-identified system of equations, whether the model is correctly specified or not. Hence, in this formulation the standard asymptotic theory of just-identified GMM estimators applies [46] (see also [32] for the application of this idea to ET under misspecification). The asymptotic variance of a just-identified GMM of the form $n^{-1} \sum_i \phi(x_i, \hat{\beta}) = 0$ is given by

$$(26) \qquad (E[\partial \phi'(x_i, \beta)/\partial \beta])^{-1} (E[\phi(\beta)\phi'(\beta)])(E[\partial \phi'(x_i, \beta)/\partial \beta])^{-1}.$$

For $k$ sufficiently large, we can always choose $F_k(x)$ so as to satisfy the necessary regularity conditions for this expression to hold. In particular, the compact support of $F_k(x)$ enables $E[g(x_i, \theta)/(1 - \lambda' g(x_i, \theta))]$ to exist for $(\theta', \lambda')'$ in some neighborhood of the pseudo-true value $(\theta_k^{*\prime}, \lambda_k^{*\prime})'$. The asymptotic distribution of $(\hat{\theta}'_{k,n}, \hat{\lambda}'_{k,n})'$ is then given by

$$n^{1/2}((\hat{\theta}'_{k,n}, \hat{\lambda}'_{k,n}) - (\theta_{k,n}^{*\prime}, \lambda_{k,n}^{*\prime}))' \xrightarrow{d} N(0, H_k^{-1} S_k H_k^{-1}), \tag{27}$$

$$\text{as } n \to \infty \text{ for fixed } k,$$

where

$$S_k = E\left[\begin{pmatrix} \tau_i^2 G_i' \lambda_k^* \lambda_k^{*\prime} G_i & \tau_i^2 G_i' \lambda_k^* g_i' \\ \tau_i^2 g_i \lambda_k^{*\prime} G_i & \tau_i^2 g_i g_i' \end{pmatrix}\right], \tag{28}$$

$$H_k = E\left[\begin{pmatrix} \dot{\tau}_i G_i' \lambda_k^* \lambda_k^{*\prime} G_i + \tau_i \frac{\partial (G_i' \lambda_k^{*\prime})}{\partial \theta'} & \dot{\tau}_i G_i' \lambda_k^* g_i' + \tau_i G_i' \\ \dot{\tau}_i g_i \lambda_k^{*\prime} G_i + \tau_i G_i & \dot{\tau}_i g_i g_i' \end{pmatrix}\right], \tag{29}$$

and where $\tau_i = \tau(\lambda_k^{*\prime} g_i) = (1 - \lambda_k^{*\prime} g_i)^{-1}$, $\dot{\tau}_i = \frac{\partial \tau(\xi)}{\partial \xi}|_{\xi = \lambda' g_i} = (1 - \lambda_k^{*\prime} g_i)^{-2} = \tau_i^2$ and where all moments are evaluated at $\theta_k^*$ and $\lambda_k^* = \text{plim}_{n \to \infty} \hat{\lambda}_{k,n}$ and assuming that $x$ is drawn from $F_k(x)$ (i.e., $E[\cdot] \equiv E_{F_k}[\cdot]$).

We focus on the upper left $N_\theta \times N_\theta$ submatrix of $H_k^{-1} S_k H_k^{-1}$, denoted by $\Sigma_k$. For a given $k$, the submatrix $\Sigma_k$ provides the asymptotic variance of $\hat{\theta}_{k,n}$. We will now analyze the behavior of $\Sigma_k$ as $k \to \infty$ (we are not claiming that this provides the asymptotic variance of EL for infinite support at this point). Since EL's implied probabilities must be positive (see, e.g., [4, 50]), it follows that $(1 - \lambda_k^{*\prime} g(x, \theta_k^*))^{-1} > 0$ for all $x \in \mathcal{C}_k$, or

$$\max_{x \in \mathcal{C}_k}(\lambda_k^{*\prime} g(x, \theta_k^*)) < 1. \tag{30}$$

Since $\{g(x, \theta_k^*) : x \in \mathcal{X}\}$ is unbounded in every direction, the set $\{g(x, \theta_k^*) : \in \mathcal{C}_k\}$ becomes unbounded in every direction as $k \to \infty$. Hence, the only way to satisfy equation (30) is to have $\lambda_k^* \to 0$ as $k \to \infty$. Since $\lambda_k^* \to 0$ as $k \to \infty$, the expressions for $S_k$ and $H_k$ can be simplified by noting that when the product $H_k^{-1} S_k H_k^{-1}$ is calculated, any term containing $\lambda_k^*$ will be dominated by terms not containing $\lambda_k^*$. We then obtain [keeping the $\tau_i = \tau(\lambda_k^{*\prime} g_i)$ prefactors even though $\lambda_k^* \to 0$ because the $g(x, \theta_k^*)$ are unbounded and it is not clear whether we necessarily have $\tau_i \to 1$]

$$S_k \to \begin{bmatrix} 0 & 0 \\ 0 & E[\tau_i^2 g_i g_i'] \end{bmatrix},$$

$$H_k^{-1} \to \begin{bmatrix} 0 & E[\tau_i G_i'] \\ E[\tau_i G_i] & E[\tau_i^2 g_i g_i'] \end{bmatrix}^{-1} \equiv \begin{bmatrix} B_{11} & B_{12} \\ B_{21} & B_{22} \end{bmatrix}. \tag{31}$$

(Note that the sequence $F_k$ can be easily chosen so that the smallest eigenvalue $H_k$ remains bounded away from zero for all $k$ sufficiently large, since the moment conditions remain the same over $k$ and $\mathcal{G}_k$ increases with $k$. Hence, $\lim_{k\to\infty} H_k$ can be assumed nonsingular and interchanging the limit as $k \to \infty$ and the matrix inversion operation is justified.) We then have that

$$\Sigma_k = B_{12} E[\tau_i^2 g_i g_i'] B_{21} + \rho_k, \tag{32}$$

where $\rho_k$ is a remainder that vanishes as $k \to \infty$ (its precise form has no bearing on the rest of the argument). By the partitioned inverse formula,

$$B_{21} = (E[\tau_i^2 g_i g_i'])^{-1} E[\tau_i G_i] (E[\tau_i G_i'] (E[\tau_i^2 g_i g_i'])^{-1} E[\tau_i G_i])^{-1} = B_{12}'. \tag{33}$$

Substituting this expression for $B_{21}$ into equation (34) yields

$$\Sigma_k = (E[\tau_i G_i'] (E[\tau_i^2 g_i g_i'])^{-1} E[\tau_i G_i])^{-1} + \rho_k. \tag{34}$$

We will now show that $\Sigma_k$ diverges as $k \to \infty$. For EL, $\lambda_k^*$ is such that $E[g(x_i, \theta_k^*)/(1 - \lambda_k^{*\prime} g(x_i, \theta_k^*))] = 0$. Since $E[g(x_i, \theta_k^*)/(1 - \lambda_k^{*\prime} g(x_i, \theta_k^*))] = E[g(x_i, \theta_k^*)] + E[g(x_i, \theta_k^*)g'(x_i, \theta_k^*)/(1 - \lambda_k^{*\prime} g(x_i, \theta_k^*))\lambda_k^*]$, we have

$$\Omega_k \lambda_k^* = -E[g(x_i, \theta_k^*)], \tag{35}$$

where $\Omega_k = E[g(x_i, \theta_k^*) g'(x_i, \theta_k^*)/(1 - \lambda_k^{*\prime} g(x_i, \theta_k^*))]$. Since $\inf_{k \geq \bar{k}} E[g(x_i, \theta_k^*)] > 0$ for some $\bar{k} \in \mathbb{N}$ by construction, having $\lambda_k^* \to 0$ as $k \to \infty$ is only possible if at least one of the eigenvalues of $\Omega_k$ diverges as $k \to \infty$. Let $v$ be a (unit) eigenvector associated with one of these eigenvalues. Then, by the CSI $v' \Omega_k v$ equals

$$E\left[\frac{v'g(x_i, \theta_k^*)}{(1 - \lambda_k^{*\prime} g(x_i, \theta_k^*))} v' g(x_i, \theta_k^*)\right]$$
$$\leq \left(E\left[\frac{(v'g(x_i, \theta_k^*))^2}{(1 - \lambda_k^{*\prime} g(x_i, \theta_k^*))^2}\right] E[(v' g(x_i, \theta_k^*))^2]\right)^{1/2}. \tag{36}$$

Since $E[(v'g(x_i, \theta_k^*))^2] \leq \sup_{\theta \in \Theta} E[\|g(x_i, \theta)\|^2] < \infty$, (36) therefore implies that $E[\frac{(v'g(x_i, \theta_k^*))^2}{(1-\lambda_k^{*\prime} g(x_i, \theta_k^*))^2}] = E[\tau_i^2 v' g_i g_i' v]$ diverges and thus that $E[\tau_i^2 g_i g_i']$ has a divergent eigenvalue. Since $E[\tau_i^2 g_i g_i']$ enters the expression of $\Sigma_k$ [given by equation (34)], $\Sigma_k$ has at least one divergent eigenvalue as $k \to \infty$. Note that the other terms entering the expression of $\Sigma_k$ cannot compensate for the explosive behavior of $E[\tau_i^2 g_i g_i']$, since a simple application of the CSI shows that, as $k \to \infty$, $\|E[\tau_i G_i]\| = \|E[(1 + \tau_i \lambda_k^{*\prime} g_i) G_i]\| = O(E[\tau_i \|g_i\| \|\lambda_k^*\| \|G_i\|]) = O((E[\tau_i^2 \|g_i\|^2])^{1/2}) (E[\|G_i\|^2])^{1/2} \|\lambda_k^*\| = o((E[\tau_i^2 \|g_i\|^2])^{1/2}) = o((E[\tau_i^2 \times \|g_i g_i'\|])^{1/2}) = o((E[\tau_i^2 v' g_i g_i' v])^{1/2})$.

We will now show that the divergent behavior of $\Sigma_k$ implies that EL is not root $n$ consistent. We start by calculating the probability that $\hat{\theta}_{k,n}$ lies outside of a root $n$ neighborhood of the pseudo-true value $\theta_k^*$. Let $P_{k,n}$ be the finite sample distribution of $n^{1/2}(u'\Sigma_k u)^{-1/2} u'(\hat{\theta}_{k,n} - \theta_k^*)$ for some conformable unit vector $u$ such that $u'\Sigma_k u \to \infty$ as $k \to \infty$ ($u$ is an eigenvector associated with one of the divergent eigenvalues of $\Sigma_k$). Let $P_{k,\infty}$ denote the corresponding asymptotic distribution, the c.d.f. of a $N(0,1)$ for all $k$. For a given $\xi < 0$, the probability that $u'(\hat{\theta}_{k,n} - \theta_k^*) \leq n^{-1/2}\xi$ is $P_{k,n}((u'\Sigma_k u)^{-1/2}\xi)$.

Let $n_k = \min\{n : \sup_{m \geq n} |P_{k,m}((u'\Sigma_k u)^{-1/2}\xi) - P_{k,\infty}((u'\Sigma_k u)^{-1/2}\xi)| \leq k^{-1}\}$. This defines the sample size beyond which the difference [at $(u'\Sigma_k u)^{-1/2}\xi$] between the finite sample and asymptotic distribution is less than $k^{-1}$. Such a finite $n$ can always be found, since $P_{k,n}$ converges pointwise to $P_{k,\infty}$. Now define the "inverse" sequence $k_n = \max\{k : n_k \leq n\}$. Note that $k_n \to \infty$, as $n \to \infty$, since $n_k \to \infty$ as $k \to \infty$.

Since $P[|u'(\hat{\theta}_{\infty,n} - \theta_\infty^*)| \geq \varepsilon] \geq P[|u'(\hat{\theta}_{k,n} - \theta_k^*)| \geq \varepsilon]$ for all $\varepsilon > 0$ and any $n$, by the construction of $F_k$, $P[u'(\hat{\theta}_{\infty,n} - \theta_\infty^*) \leq n^{-1/2}\xi] \geq P[u'(\hat{\theta}_{k,n} - \theta_k^*) \leq n^{-1/2}\xi] = P_{k,n}((u'\Sigma_k u)^{-1/2}\xi)$ for any $k$ and any $n$ and any $\xi < 0$. In particular, for $k = k_n$,

$$\begin{align}
P[u'(\hat{\theta}_{\infty,n} - \theta_\infty^*) &\leq n^{-1/2}\xi] \geq P_{k_n,n}((u'\Sigma_{k_n} u)^{-1/2}\xi) \tag{37}\\
&= P_{k_n,\infty}((u'\Sigma_{k_n} u)^{-1/2}\xi) \\
&\quad + (P_{k_n,n}((u'\Sigma_{k_n} u)^{-1/2}\xi) - P_{k_n,\infty}((u'\Sigma_{k_n} u)^{-1/2}\xi)) \tag{38}\\
&\geq P_{k_n,\infty}((u'\Sigma_{k_n} u)^{-1/2}\xi) - k_n^{-1} \tag{39}
\end{align}$$

by the definition of $k_n$. As $n \to \infty$, $k_n^{-1} \to 0$. Since $P_{k,\infty}$ is the same for all $k$ and is continuous [it is the c.d.f. of a $N(0,1)$], for any $k$ we have $\lim_{n \to \infty} P_{k_n,\infty}((u'\Sigma_{k_n} u)^{-1/2}\xi) = \lim_{n \to \infty} P_{k,\infty}((u'\Sigma_{k_n} u)^{-1/2}\xi) = P_{k,\infty} \times (\lim_{n \to \infty}(u'\Sigma_{k_n} u)^{-1/2}\xi) = P_{k,\infty}(0) = 1/2$, where we have used the fact that $(u'\Sigma_{k_n} u)^{-1/2} \to 0$ since $u'\Sigma_k u$ diverges as $k \to \infty$. We then have $\lim_{n \to \infty} P[u'(\hat{\theta}_{\infty,n} - \theta_\infty^*) \leq n^{-1/2}\xi] \geq 1/2$ for any $\xi < 0$. A similar reasoning for $\xi > 0$ implies that $\lim_{n \to \infty} P[u'(\hat{\theta}_{\infty,n} - \theta_\infty^*) \geq n^{-1/2}\xi] \geq 1/2$. It follows that $\hat{\theta}_{\infty,n}$ lies outside a $n^{-1/2}$ neighborhood of $\theta_\infty^*$ with probability approaching $1/2 + 1/2 = 1$ as $n \to \infty$, thus ruling out root $n$ convergence.

To summarize, for any EL estimator $\hat{\theta}_{\infty,n}$ based on a distribution $F_\infty(x)$ with unbounded support, there exists a family of other estimators $\hat{\theta}_{k,n}$ based on compactly supported distributions $F_k(x)$ all having a narrower distribution than EL for each $n$. Yet the asymptotic variance of $\hat{\theta}_{k,n}$ diverges as $k \to \infty$. By a standard diagonal argument, there exists an estimator sequence $\hat{\theta}_{k_n,n}$ that is not root $n$ consistent but whose distribution is narrower than the one of EL at each $n$. Hence EL is not root $n$ consistent. $\square$

PROOF OF THEOREM 2.

$$\ln \hat{L} \equiv n^{-1} \sum_i \ln n\hat{w}_i = n^{-1} \sum_i \ln\left( \exp(\hat{\lambda}' g_i) \Big/ \left(n^{-1} \sum_j \exp(\hat{\lambda}' g_j)\right) \right)$$

$$= n^{-1} \sum_i \hat{\lambda}' g_i - \ln\left(n^{-1} \sum_j \exp(\hat{\lambda}' g_j)\right)$$

$$= -\ln\left(n^{-1} \sum_j \exp(\hat{\lambda}'(g_j - \hat{g}))\right),$$

$$\frac{d\ln \hat{L}}{d\theta'} = n^{-1} \sum_i \frac{d(\hat{\lambda}' g_i)}{d\theta'} - \left(n^{-1} \sum_j \exp(\hat{\lambda}' g_j)\right)^{-1} n^{-1} \sum_i \exp(\hat{\lambda}' g_i) \frac{d(\hat{\lambda}' g_i)}{d\theta'}$$

$$= n^{-1} \sum_i \frac{d(\hat{\lambda}' g_i)}{d\theta'} - n^{-1} \sum_i n\hat{w}_i \frac{d(\hat{\lambda}' g_i)}{d\theta'} = n^{-1} \sum_i (1 - n\hat{w}_i) \frac{d(\hat{\lambda}' g_i)}{d\theta'}.$$

From (15), the first order condition for $\hat{\lambda}$ is $\sum_i g_i \exp(\hat{\lambda}' g_i) = 0$. $\square$

PROOF OF THEOREM 3. Expanding the ETEL first-order conditions for $\hat{\theta}$ and $\hat{\lambda}$ around $\theta = \theta^*$ and $\lambda = 0$ reveals an expansion identical to that of EL at least up to $O_p(n^{-1/2})$ in an $O(n^{-1/2})$ neighborhood of $\theta = \theta^*$ and $\lambda = 0$,

$$n^{-1} \sum_i \begin{bmatrix} 0 \\ g(x_i, \theta^*) \end{bmatrix}$$

$$+ n^{-1} \sum_i \begin{pmatrix} 0 & G'(x_i, \theta^*) \\ G(x_i, \theta^*) & g(x_i, \theta^*)g'(x_i, \theta^*) \end{pmatrix} \begin{bmatrix} \theta - \theta^* \\ \lambda \end{bmatrix} = o_p(n^{-1/2}).$$

Calculational details can be found in [57]. In addition, in a $O(n^{-1/2})$ neighborhood of $\theta = \theta^*$, both the ETEL and the EL objective functions for $\hat{\theta}$ share the same expansion in $(\theta - \theta^*)$ at least up to $O_p(n^{-1})$,

$$-\tfrac{1}{2}(\theta - \theta^*)' \left(n^{-1} \sum_i G'(x_i, \theta^*)\right) (\hat{\Omega}^*)^{-1} \left(n^{-1} \sum_i G(x_i, \theta^*)\right) (\theta - \theta^*) + o_p(n^{-1}),$$

where $\hat{\Omega}^* = n^{-1} \sum_i g(x_i, \theta^*) g'(x_i, \theta^*)$. It is known (see, e.g., [47, 49]) that the EL estimator is asymptotically such that the solutions $\hat{\lambda}_{\text{EL}}$ and $\hat{\theta}_{\text{EL}}$ lie within the $O(n^{-1/2})$ neighborhood where the remainder terms of these expansions are negligible. Hence, asymptotically, $\hat{\lambda}_{\text{EL}}$ and $\hat{\theta}_{\text{EL}}$ also solve the ETEL first-order conditions, apart from negligible remainders. Since $\hat{\lambda}_{\text{EL}} \xrightarrow{p} 0$, and since both the EL and the ETEL objective functions for $\hat{\theta}$ converge to their maximum possible value when $\hat{\lambda}_{\text{EL}} \xrightarrow{p} 0$ and $\hat{\lambda}_{\text{ETEL}} \xrightarrow{p} 0$, respectively, the existence of another solution outside of the neighborhood of validity of the above expansions can be ruled out. ETEL thus inherits all the first-order properties of EL established in [47, 49]. □

PROOF OF THEOREM 4. The first conclusion follows from the fact that the implied probabilities are given by

$$\hat{w}_i(\theta) = \exp(\hat{\lambda}(\theta)'g(x_i,\theta)) \Big/ \left(\sum_j \exp(\hat{\lambda}(\theta)'g(x_j,\theta))\right),$$

a necessarily positive quantity for any $\hat{\lambda}$ and $\theta$. The second conclusion holds for any estimator where $\theta$ is the extremum of a differentiable objective function:

$$\frac{\partial \ln \hat{L}(T(\beta))}{\partial \beta} = \frac{\partial T(\beta)'}{\partial \beta} \frac{\partial \ln \hat{L}(\theta)}{\partial \theta} = 0$$

if and only if $\partial \ln \hat{L}(\theta)/\partial \theta = 0$ since $\partial T(\beta)'/\partial \beta$ has full rank [$T(\beta)$ being one-to-one]. The third conclusion can be shown by noting that any invertible linear transformation of the moment function $A(\theta)g(x_i,\theta)$ simply causes the Lagrange multiplier $\hat{\lambda}(\theta)$ to become $((A(\theta))^{-1})'\lambda(\theta)$. Indeed, under these two transformations, the first-order conditions for both $\hat{\theta}$ and $\hat{\lambda}(\hat{\theta})$ remain satisfied,

$$n^{-1}\sum_i (1 - n\hat{w}_i(\hat{\theta}))\,d(\hat{\lambda}'(A(\theta))^{-1}A(\theta)g_i)/d\theta$$

$$= n^{-1}\sum_i (1 - n\hat{w}_i(\hat{\theta}))\,d(\hat{\lambda}'g_i)/d\theta = 0,$$

where $\hat{w}_i(\hat{\theta}) = \exp(\hat{\lambda}'(A(\theta))^{-1}A(\theta)g_i)/(\sum_j \exp(\hat{\lambda}'(A(\theta))^{-1}A(\theta)g_j)) = \exp(\hat{\lambda}'g_i)/(\sum_j \exp(\hat{\lambda}'g_j))$ and $n^{-1}\sum_i \exp(\hat{\lambda}'(A(\theta))^{-1}A(\theta)g_i)A(\theta)g_i = A(\theta)n^{-1}\sum_i \exp(\hat{\lambda}'g_i)g_i = 0$ if and only if $n^{-1}\sum_i \exp(\hat{\lambda}'g_i)g_i = 0$ since $A(\theta)$ is invertible. □

PROOF OF THEOREM 5. By Theorem 2, the first order condition for $\hat{\theta}$ is

$$\frac{d\ln \hat{L}}{d\theta'} = n^{-1}\sum_i (1-n\hat{w}_i)\left(g_i'\frac{\partial \hat{\lambda}}{\partial \theta'} + \hat{\lambda}'G_i\right)$$

(40)
$$= n^{-1}\sum_i g_i'\frac{\partial \hat{\lambda}}{\partial \theta'} + \hat{\lambda}'n^{-1}\sum_i G_i - \left(\sum_i \hat{w}_i g_i'\right)\frac{\partial \hat{\lambda}}{\partial \theta'} - \hat{\lambda}'\sum_i \hat{w}_i G_i$$

$$= \hat{g}'\frac{\partial \hat{\lambda}}{\partial \theta'} + \hat{\lambda}'\hat{G} - 0 - \hat{\lambda}'\tilde{G}.$$

To find $\partial \hat{\lambda}/\partial \theta'$, we note that the total differential of $\sum_i \exp(\hat{\lambda}' g_i) g_i = 0$ yields

$$\sum_i \exp(\hat{\lambda}' g_i) g_i g_i' \, d\hat{\lambda} + \sum_i \exp(\hat{\lambda}' g_i) G_i \, d\theta + \sum_i g_i \exp(\hat{\lambda}' g_i) \hat{\lambda}' G_i \, d\theta = 0,$$

$$\sum_i g_i g_i' \hat{w}_i \, d\hat{\lambda} + \sum_i \hat{w}_i (I + g_i \hat{\lambda}') G_i \, d\theta = 0,$$

implying that

$$(41) \qquad \frac{\partial \hat{\lambda}}{\partial \theta'} = -\tilde{\Omega}^{-1} \left( \sum_i \hat{w}_i (I + g_i \hat{\lambda}') G_i \right).$$

Substituting this result into equation (40) gives

$$\frac{\partial \ln \hat{L}(\theta)}{\partial \theta'} = -\hat{g}' \tilde{\Omega}^{-1} \left( \sum_i \hat{w}_i (I + g_i \hat{\lambda}') G_i \right) + \hat{\lambda}' \hat{G} - \hat{\lambda}' \tilde{G}$$

$$(42) \qquad = -\hat{g}' \tilde{\Omega}^{-1} \tilde{G} - \hat{g}' \tilde{\Omega}^{-1} \sum_i \hat{w}_i g_i \hat{\lambda}' G_i + \hat{\lambda}' \hat{G} - \hat{\lambda}' \tilde{G}$$

$$= -\hat{g}' \tilde{\Omega}^{-1} \tilde{G} - \hat{g}' \tilde{\Omega}^{-1} \sum_i \hat{w}_i g_i \hat{\lambda}' G_i + n^{-1} \sum_i (1 - n\hat{w}_i) \hat{\lambda}' G_i.$$

By the first-order equivalence between EL and ETEL established in Theorem 3 and using Theorem 3.1 in [47], $\hat{\lambda}(\hat{\theta}) = O_p(n^{-1/2})$ and $\hat{g} = O_p(n^{-1/2})$ for $\hat{\theta}$ such that $\|\theta - \theta^*\| = O_p(n^{-1/2})$. These facts, along with the fact that $\sup_{\theta \in \Theta} \max_{i \leq n} \|g_i\| = o_p(n^{1/2})$ by part 4 of Assumption 1, provide us with asymptotic expansions for $n\hat{w}_i$ and $\hat{\lambda}$,

$$n\hat{w}_i = \frac{\exp(\hat{\lambda}' g_i)}{n^{-1} \sum_j \exp(\hat{\lambda}' g_j)}$$

$$(43) \qquad = \frac{1 + \hat{\lambda}' g_i + O((\hat{\lambda}' g_i)^2)}{1 + \hat{\lambda}' \hat{g} + O_p(n^{-1})}$$

$$= \frac{1 + \hat{\lambda}' g_i + O_p(n^{-1}) \|g_i\|^2}{1 + O_p(n^{-1}) + O_p(n^{-1})}$$

$$= 1 + \hat{\lambda}' g_i + O_p(n^{-1}) \|g_i\|^2.$$

An expansion for $\hat{\lambda}$ is obtained by noting that the left-hand side of $n^{-1} \sum_i g_i \times \exp(g_i' \hat{\lambda}) = 0$ can be written as

$$n^{-1} \sum_i g_i (1 + g_i' \hat{\lambda}) + R_0 = n^{-1} \sum_i g_i + \left( n^{-1} \sum_i g_i g_i' \right) \hat{\lambda} + R_0$$

$$= n^{-1} \sum_i g_i + \left( n^{-1} \sum_i n\hat{w}_i g_i g_i' \right) \hat{\lambda} + R_0 + R_1$$

$$= \hat{g} + \tilde{\Omega}\hat{\lambda} + R_0 + R_1,$$

implying that

(44) $$\hat{\lambda} = -\tilde{\Omega}^{-1}\hat{g} - \tilde{\Omega}^{-1}(R_0 + R_1),$$

where the remainder terms $R_0, R_1$ can be bounded using the assumption $E[\sup_{\theta \in \mathcal{N}} \|g_i\|^4] < \infty$ and (43): $\|R_0\| = O_p(n^{-1})n^{-1}\sum_i \|g_i\|^3 = O_p(n^{-1})$ and $\|R_1\| \leq n^{-1}\sum_i (n\hat{w}_i - 1)\|g_i\|^2 \|\hat{\lambda}\| = n^{-1}\sum_i O(\|\hat{\lambda}\|\|g_i\|)\|g_i\|^2 \times \|\hat{\lambda}\| = O(\|\hat{\lambda}\|^2)n^{-1}\sum_i \|g_i\|^3 = O_p(n^{-1})$.

Substituting the expansion (43) into the last term of (42) yields

(45) 
$$\frac{\partial \ln \hat{L}(\theta)}{\partial \theta'} = -\hat{g}'\tilde{\Omega}^{-1}\tilde{G} - \hat{g}'\tilde{\Omega}^{-1}\sum_i \hat{w}_i g_i \hat{\lambda}' G_i$$
$$+ n^{-1}\sum_i \hat{\lambda}' g_i \hat{\lambda}' G_i + R_2,$$

where $\|R_2\| \leq O_p(n^{-1})n^{-1}\sum_i \|g_i\|^2 \|\hat{\lambda}\|\|G_i\| \leq O_p(n^{-3/2})n^{-1}\sum_i \|g_i\|^2\|G_i\| \leq O_p(n^{-3/2})(n^{-1}\sum_i \|g_i\|^4)^{1/2}(n^{-1}\sum_i \|G_i\|^2)^{1/2} = O_p(n^{-3/2})$, after using the CSI and the facts that $E[\sup_{\theta \in \mathcal{N}} \|g_i\|^4] < \infty$ and $E[\sup_{\theta \in \mathcal{N}} \|G_i\|^2] < \infty$. Then (45) becomes

$$\frac{\partial \ln \hat{L}(\theta)}{\partial \theta'} = -\hat{g}'\tilde{\Omega}^{-1}\tilde{G} - \hat{g}'\tilde{\Omega}^{-1}\sum_j \hat{w}_j g_j \hat{\lambda}' G_j$$
$$- (\hat{\lambda})' n^{-1}\sum_j g_j \hat{\lambda}' G_j + O_p(n^{-3/2}),$$

where the $\hat{\lambda}$ in parentheses can be replaced by expansion (44),

(46)
$$\frac{\partial \ln \hat{L}(\theta)}{\partial \theta'} = -\hat{g}'\tilde{\Omega}^{-1}\tilde{G} - \hat{g}'\tilde{\Omega}^{-1}\sum_j \hat{w}_j g_j \hat{\lambda}' G_j$$
$$+ \hat{g}'\tilde{\Omega}^{-1} n^{-1}\sum_j g_j \hat{\lambda}' G_j + R_3 + O_p(n^{-3/2}),$$

where $\|R_3\| = O_p(n^{-1})n^{-1}\sum_j \|g_j\|\|\hat{\lambda}\|\|G_j\| = O_p(n^{-3/2})n^{-1}\sum_j \|g_j\|\|G_j\| \leq O_p(n^{-3/2})(n^{-1}\sum_j \|g_j\|^2)^{1/2}(n^{-1}\sum_j \|G_j\|^2)^{1/2} = O_p(n^{-3/2})$ by the CSI, $E[\sup_{\theta \in \mathcal{N}} \|g_i\|^4] < \infty$ and $E[\sup_{\theta \in \mathcal{N}} \|G_i\|^2] < \infty$. Then (46) becomes

$$\frac{\partial \ln \hat{L}(\theta)}{\partial \theta'} = -\hat{g}'\tilde{\Omega}^{-1}\tilde{G} + \hat{g}'\tilde{\Omega}^{-1} n^{-1}\sum_j (1 - n\hat{w}_j) g_j \hat{\lambda}' G_j$$

(47)
$$+ O_p(n^{-3/2})$$
$$= -\hat{g}'\tilde{\Omega}^{-1}\tilde{G} - \hat{g}'\tilde{\Omega}^{-1}n^{-1}\sum_j(\hat{\lambda}'g_j)g_j\hat{\lambda}'G_j + R_4$$
$$+ O_p(n^{-3/2}),$$

where we have used the expansion (43) again and where $\|R_4\| \leq \|\hat{g}'\tilde{\Omega}^{-1}\|n^{-1} \times \sum_j O((\hat{\lambda}'g_j)^2)\|g_j\|\|\hat{\lambda}\|\|G_j\| = \|\hat{g}\|\|\hat{\lambda}\|^3\|\tilde{\Omega}^{-1}\|n^{-1}\sum_j \|g_j\|^2\|g_j\|\|G_j\| \leq O_p(n^{-2})(\max_{i\leq n}\|g_j\|)n^{-1}\sum_j \|g_j\|^2\|G_j\| = O_p(n^{-2})O_p(n^{1/2})n^{-1}\sum_j \|g_j\|^2 \times \|G_j\| = O_p(n^{-3/2})$ by the CSI, the assumptions that $E[\sup_{\theta \in \mathcal{N}} \|g_i\|^4] < \infty$ and $E[\sup_{\theta \in \mathcal{N}} \|G_i\|^2] < \infty$ and the fact that $E[\|g_i\|^2] < \infty \Rightarrow \max_{i \leq n} \|g_i\| = O_p(n^{1/2})$ (as in [47], Lemma A1). Finally (47) becomes

$$\frac{\partial \ln \hat{L}(\theta)}{\partial \theta'} = -\hat{g}'\tilde{\Omega}^{-1}\tilde{G} + O_p(n^{-3/2}).$$

Now, the term $\hat{g}'\tilde{\Omega}^{-1}\tilde{G}$ is similar to the first-order conditions for EL, except that the weights used in $\tilde{\Omega}$ and $\tilde{G}$ are the ET rather than the EL weights. However, by (43) and a similar expansion for the EL weights, $n(\hat{w}_{i,\mathrm{ET}} - \hat{w}_{i,\mathrm{EL}}) = O_p(n^{-1})\|g_i\|^2$. This fact, along with $\hat{g} = O_p(n^{-1/2})$, implies that

$$\hat{g}'\left(\sum_i \hat{w}_{i,\mathrm{ET}} g_i g_i'\right)^{-1}\left(\sum_i \hat{w}_{i,\mathrm{ET}} G_i\right)$$
$$= \hat{g}'\left(n^{-1}\sum_i n\hat{w}_{i,\mathrm{EL}} g_i g_i' + R_5\right)^{-1}\left(n^{-1}\sum_i n\hat{w}_{i,\mathrm{EL}} G_i + R_6\right)$$
$$= \hat{g}'\left(\sum_i \hat{w}_{i,\mathrm{EL}} g_i g_i'\right)^{-1}\left(\sum_i \hat{w}_{i,\mathrm{EL}} G_i\right) + O_p(n^{-1/2})O_p(n^{-1})$$

by the differentiability of the inverse and the fact that $\|R_5\| \leq n^{-1}\sum_i O_p(n^{-1}) \times \|g_i\|^2\|g_i\|^2 = O_p(n^{-1})n^{-1}\sum_i \|g_i\|^4 = O_p(n^{-1})$ and $\|R_6\| \leq n^{-1} \times \sum_i O_p(n^{-1})\|g_i\|^2\|G_i\| = O_p(n^{-1})$.

This implies that the first-order condition for ETEL is the same as that of EL up to $O_p(n^{-3/2})$. The continuous differentiability of $\hat{g}$ in $\theta$ implies $\|\hat{\theta}_{\mathrm{ETEL}} - \hat{\theta}_{\mathrm{EL}}\| = O_p(n^{-3/2})$ by a standard expansion of the first-order condition around $\theta = \theta^*$. □

PROOF OF THEOREM 7. Lemma A4 in [47] establishes that under regularity conditions implied by the ones given in the statement of the present theorem, a just-identified GMM estimator $\hat{\beta}$ defined by $n^{-1}\sum_i \phi(x_i, \hat{\beta}) = 0$ admits a stochastic expansion of the form

(48) $$\hat{\beta}_l - \beta_l^* = n^{-1/2}\bar{\Psi}_l + n^{-1}\bar{Q}_l + n^{-3/2}\bar{R}_l + O_p(n^{-2}),$$

where

$$\bar{Q}_l = \sum_j \bar{\Psi}_{l,j} \bar{\Psi}_j + \tfrac{1}{2} \sum_{j,k} \Psi_{l,jk} \bar{\Psi}_j \bar{\Psi}_k,$$

$$\bar{R}_l = \sum_j \bar{\Psi}_{l,j} \bar{Q}_j + \sum_{j,k} \Psi_{l,jk} \bar{\Psi}_j \bar{Q}_k$$

$$+ \tfrac{1}{2} \sum_{j,k} \bar{\Psi}_{l,jk} \bar{\Psi}_j \bar{\Psi}_k + \tfrac{1}{6} \sum_{j,k,h} \Psi_{l,jkh} \bar{\Psi}_j \bar{\Psi}_k \bar{\Psi}_h,$$

$$\bar{\Psi}_l = \sum_q \Phi_{lq}^{-1} \bar{\Phi}_q, \qquad \bar{\Psi}_{l,j} = \sum_q \Phi_{lq}^{-1} \bar{\Phi}_{q,j}, \qquad \bar{\Psi}_{l,jk} = \sum_q \Phi_{lq}^{-1} \bar{\Phi}_{q,jk},$$

$$\Psi_l = \sum_q \Phi_{lq}^{-1} \Phi_q, \qquad \Psi_{l,j} = \sum_q \Phi_{lq}^{-1} \Phi_{q,j}, \qquad \Psi_{l,jk} = \sum_q \Phi_{lq}^{-1} \Phi_{q,jk},$$

$$\Psi_{l,jkh} = \sum_q \Phi_{lq}^{-1} \Phi_{q,jkh},$$

$$\Phi^{-1} = \left( E\left[ \left. \frac{\partial \phi(x_i, \beta)}{\partial \beta'} \right|_{\beta=\beta^*} \right] \right)^{-1},$$

$$\Phi_{l,j} = E\left[ \left. \frac{\partial \phi_l(x_i, \beta)}{\partial \beta_j} \right|_{\beta=\beta^*} \right], \qquad \Phi_{l,jk} = \left[ \left. \frac{\partial^2 \phi_l(x_i, \beta)}{\partial \beta_j \partial \beta_k} \right|_{\beta=\beta^*} \right],$$

$$\Phi_{l,jkh} = \left[ \left. \frac{\partial^3 \phi_l(x_i, \beta)}{\partial \beta_j \partial \beta_k \partial \beta_h} \right|_{\beta=\beta^*} \right],$$

$$\bar{\Phi}_l = n^{-1/2} \sum_i \phi_l(x_i, \beta^*), \qquad \bar{\Phi}_{l,j} = n^{-1/2} \sum_i \left( \left. \frac{\partial \phi_l(x_i, \beta)}{\partial \beta_j} \right|_{\beta=\beta^*} - \Phi_{l,j} \right),$$

$$\bar{\Phi}_{l,jk} = n^{-1/2} \sum_i \left( \left. \frac{\partial^2 \phi_l(x_i, \beta)}{\partial \beta_j \partial \beta_k} \right|_{\beta=\beta^*} - \Phi_{l,jk} \right).$$

(We have adapted Newey and Smith's result to follow our notation and slightly simplified it using the fact that $\Psi_{l,jk} = \Psi_{l,kj}$.) We now write the ETEL and EL estimators as just identified GMM estimators that can be easily compared. As shown in Lemma 9, and as discussed in Section 3.2.3 in the text, the ETEL estimator can be written as a subvector $\hat{\theta}$ of an augmented parameter vector $\hat{\beta} = (\hat{\tau}, \hat{\kappa}', \hat{\lambda}', \hat{\theta}')'$ that solves a just-identified vector of moment conditions $n^{-1} \sum_i \phi^{\text{ETEL}}(x_i, \hat{\beta}_{\text{ETEL}}) = 0$, where $\phi^{\text{ETEL}}(x_i, \hat{\beta})$ is given by (23).

It is well known that EL can also be written as a subvector $\hat{\theta}$ of an augmented parameter vector $(\hat{\kappa}', \hat{\theta}')'$ that solves a just-identified vector of moment conditions

$$n^{-1} \sum_i \begin{bmatrix} \hat{\varepsilon}_i g_i \\ \hat{\varepsilon}_i G_i' \kappa \end{bmatrix} = 0, \tag{49}$$

where $\hat{\varepsilon}_i = (1 - \hat{\kappa}' g_i)^{-1}$ and $\hat{\kappa}$ is the Lagrange multiplier of the moment constraints, which has been relabelled $\hat{\kappa}$ to simplify the comparison with ETEL. Once again, to further simplify the comparison, we augment the vector in (49) by $1 + \dim \kappa$ additional moment conditions and introduce the same number of additional parameters $(\hat{\tau}, \hat{\lambda})$, where $\tau \in \mathbb{R}$ and $\lambda \in \mathbb{R}^{\dim \kappa}$,

$$n^{-1} \sum_i [(\hat{\tau}_i - \hat{\tau}) \quad \hat{\tau}_i g_i' \quad \hat{\varepsilon}_i g_i' \quad (\hat{\varepsilon}_i G_i' \hat{\kappa})']' = 0,$$

where $\hat{\tau}_i = \exp(\hat{\lambda}' g_i)$. In this fashion, the dimension of the vector of moment conditions and the number of parameters are the same in ETEL, as in EL. The additional moment conditions merely define the values of the new parameters $(\hat{\tau}, \hat{\lambda})$ and do not change the values of $(\hat{\kappa}', \hat{\theta}')'$. Indeed, whenever $(\hat{\kappa}', \hat{\theta}')'$ are such that the bottom two subvectors are zero, one can always find a value of $(\hat{\tau}, \hat{\lambda})$ that will make the top two subvectors vanish as well. (There exists $\hat{\lambda}$ such that $n^{-1} \sum_i \hat{\tau}_i g_i = 0$ w.p.a. 1. Then, we can just set $\hat{\tau} = n^{-1} \sum_i \hat{\tau}_i$.)

Finally, since just-identified GMM is invariant under linear transformations of the vector of moment conditions, the moment conditions for EL can equivalently be written as $n^{-1} \sum_i \phi^{\text{EL}}(x_i, \hat{\beta}_{\text{EL}}) = 0$, where

$$\phi^{\text{EL}}(x_i, \hat{\beta}) = \begin{bmatrix} \hat{\tau}_i - \tau \\ \hat{\tau}_i g_i \\ \hat{\varepsilon}_i g_i - \hat{\tau}_i g_i \\ \hat{\varepsilon}_i G_i' \hat{\kappa} \end{bmatrix}. \tag{50}$$

Equipped with (23) and (50), we can construct a stochastic expansion of the form (48) for each estimator. The $O(n^{-2})$ covariance between two elements of the parameter vector, $\hat{\theta}_l$ and $\hat{\theta}_m$, is given by

$$W_{lm} \equiv \text{Covar}[\bar{Q}_{l_\theta + l}, \bar{Q}_{l_\theta + m}] + \text{Covar}[\bar{R}_{l_\theta + l}, \bar{\Psi}_{l_\theta + m}] + \text{Covar}[\bar{\Psi}_{l_\theta + l}, \bar{R}_{l_\theta + m}], \tag{51}$$

where $l_\theta = 1 + 2 \dim \lambda$. The quantities associated with each estimator will be distinguished by an "ETEL" or "EL" superscript.

We provide below the sequence of equalities that need to be established in order to show, as directly as possible, that ETEL and EL have the same $O(n^{-2})$ variance. The tedious yet straightforward calculational details that prove each statement are omitted below but can be found in [57].

(1) $\bar{\Phi}_l^{\text{ETEL}} = \bar{\Phi}_l^{\text{EL}}$ and $\Phi_{l,j}^{\text{ETEL}} = \Phi_{l,j}^{\text{EL}} \equiv \Phi_{l,j} \Rightarrow \bar{\Psi}_j^{\text{ETEL}} = \bar{\Psi}_j^{\text{EL}}$.

(2) $(1) \Rightarrow \bar{Q}_l^{\text{ETEL}} - \bar{Q}_l^{\text{EL}} = \sum_j (\bar{\Psi}_{l,j}^{\text{ETEL}} - \bar{\Psi}_{l,j}^{\text{EL}}) \bar{\Psi}_j + \frac{1}{2} \sum_{j,k} (\Psi_{l,jk}^{\text{ETEL}} - \Psi_{l,jk}^{\text{EL}}) \times \bar{\Psi}_j \bar{\Psi}_k$.

(2a) $(\bar{\Psi}_{l,j}^{\text{ETEL}} - \bar{\Psi}_{l,j}^{\text{EL}}) \bar{\Psi}_j = \sum_{q,j} \Phi_{lq}^{-1} (\bar{\Phi}_{q,j}^{\text{ETEL}} - \bar{\Phi}_{q,j}^{\text{EL}}) \bar{\Psi}_j$, where $\sum_j (\bar{\Phi}_{q,j}^{\text{ETEL}} - \bar{\Phi}_{q,j}^{\text{EL}}) \bar{\Psi}_j = 0$.

(2b) $(\Psi_{l,jk}^{\text{ETEL}} - \Psi_{l,jk}^{\text{EL}}) \bar{\Psi}_j \bar{\Psi}_k = \sum_{q,j,k} \Phi_{lq}^{-1} (\Phi_{q,jk}^{\text{ETEL}} - \Phi_{q,jk}^{\text{EL}}) \bar{\Psi}_j \bar{\Psi}_k$, where $\sum_{j,k} (\Phi_{q,jk}^{\text{ETEL}} - \Phi_{q,jk}^{\text{EL}}) \bar{\Psi}_j \bar{\Psi}_k = 0$.

(3) (2), (2a) and (2b) $\Rightarrow \bar{Q}_l^{\text{ETEL}} - \bar{Q}_l^{\text{EL}} = 0$.

(4) (1) and (3) $\Rightarrow W_{lm}^{\text{ETEL}} - W_{lm}^{\text{EL}} = \text{Covar}[\bar{R}_{l_\theta+l}^{\text{ETEL}} - \bar{R}_{l_\theta+l}^{\text{EL}}, \bar{\Psi}_{l_\theta+m}] + \text{Covar}[\bar{\Psi}_{l_\theta+l}, \bar{R}_{l_\theta+m}^{\text{ETEL}} - \bar{R}_{l_\theta+m}^{\text{EL}}]$.

(5) (1) and (3) $\Rightarrow \bar{R}_l^{\text{ETEL}} - \bar{R}_l^{\text{EL}} = \sum_j (\bar{\Psi}_{l,j}^{\text{ETEL}} - \bar{\Psi}_{l,j}^{\text{EL}}) \bar{Q}_{,j} + \sum_{j,k} (\Psi_{l,jk}^{\text{ETEL}} - \Psi_{l,jk}^{\text{EL}}) \bar{\Psi}_j \bar{Q}_k + \frac{1}{2} \sum_{j,k} (\bar{\Psi}_{l,jk}^{\text{ETEL}} - \bar{\Psi}_{l,jk}^{\text{EL}}) \bar{\Psi}_j \bar{\Psi}_k + \frac{1}{6} \sum_{j,k,h} (\Psi_{l,jkh}^{\text{ETEL}} - \Psi_{l,jkh}^{\text{EL}}) \bar{\Psi}_j \bar{\Psi}_k \bar{\Psi}_h$.

(5a) $\sum_j (\bar{\Psi}_{l_\theta+l,j}^{\text{ETEL}} - \bar{\Psi}_{l_\theta+l,j}^{\text{EL}}) \bar{Q}_{,j} = \frac{1}{2} \sum_j H_{lj} \bar{g}_j \bar{g}' P \bar{g}$, where $\bar{g} = n^{-1/2} \sum_i g_i$, $H = (G'\Omega^{-1}G)^{-1} G'\Omega^{-1}$ and $P = \Omega^{-1} - \Omega^{-1} G (G'\Omega^{-1}G)^{-1} G'\Omega^{-1}$.

(5b) $(\Psi_{l_\theta+l,jk}^{\text{ETEL}} - \Psi_{l_\theta+l,jk}^{\text{EL}}) \bar{\Psi}_k \bar{Q}_j = -\frac{1}{2} \sum_j H_{lj} \bar{g}_j \bar{g}' P \bar{g} + \Xi_{1,l}$ with $E[\Xi_l \bar{\Psi}_{l_\theta+m}] = o(n^{-2})$.

(5c) $\bar{\Psi}_j (\bar{\Psi}_{l,jk}^{\text{ETEL}} - \bar{\Psi}_{l,jk}^{\text{EL}}) \bar{\Psi}_k = 0$

(5d) $E[(\Psi_{l_\theta+l,jkh}^{\text{ETEL}} - \Psi_{l_\theta+l,jkh}^{\text{EL}}) \bar{\Psi}_j \bar{\Psi}_k \bar{\Psi}_h \bar{\Psi}_{l_\theta+m}^{\text{EL}}] = o(n^{-2})$.

(6) (5a) through (5d) $\Rightarrow \text{Covar}[\bar{R}_{l_\theta+l}^{\text{ETEL}} - \bar{R}_{l_\theta+l}^{\text{EL}}, \bar{\Psi}_{l_\theta+m}] = o(n^{-2})$.

(7) (3) and (6) $\Rightarrow$ ETEL and EL share the same $O(n^{-2})$ variance. $\square$

PROOF OF THEOREM 8. Let $g_{i,a}$ and $g_{i,b}$ denote the subvectors of $g_i$ that are mutually independent and let $\lambda_a$ and $\lambda_b$ denote the corresponding subvectors of the Lagrange multiplier. Independence holds if and only if for any measurable functions $a(g_{i,a})$ and $b(g_{i,b})$, $E[a(g_{i,a})b(g_{i,b})] = E[a(g_{i,a})]E[b(g_{i,b})]$ whenever these expectations are defined. The exponentially tilted empirical distribution estimates the moment $E[a(g_{i,a})]$ by

$$\hat{Q}_a = \left(n^{-1} \sum_j \exp(\hat{\lambda}' g_j)\right)^{-1} n^{-1} \sum_i a(g_{i,a}) \exp(\hat{\lambda}' g_i)$$

$$\xrightarrow{p} (E[\exp(\lambda' g_j)])^{-1} E[a(g_{i,a}) \exp(\lambda' g_i)]$$

$$= (E[\exp(\lambda_a' g_{i,a}) \exp(\lambda_b' g_{i,b})])^{-1} E[a(g_{i,a}) \exp(\lambda_a' g_{i,a}) \exp(\lambda_b' g_{i,b})]$$

$$= \frac{E[a(g_{i,a}) \exp(\lambda_a' g_{i,a})] E[\exp(\lambda_b' g_{i,b})]}{E[\exp(\lambda_a' g_{i,a})] E[\exp(\lambda_b' g_{i,b})]}$$

$$= \frac{E[a(g_{i,a}) \exp(\lambda_a' g_{i,a})]}{E[\exp(\lambda_a' g_{i,a})]} \equiv Q_a,$$

and similarly for $E[b(g_{i,b})]$. The exponentially tilted empirical distribution estimates the moment $E[a(g_{i,a})b(g_{i,b})]$ by

$$\hat{Q}_{ab} = \left(n^{-1}\sum_j \exp(\hat{\lambda}'g_j)\right)^{-1} n^{-1}\sum_i a(g_{i,a})b(g_{i,b})\exp(\hat{\lambda}'g_i)$$

$$\xrightarrow{p} (E[\exp(\lambda_a'g_{i,a})\exp(\lambda_b'g_{i,b})])^{-1} E[a(g_{i,a})\exp(\lambda_a'g_{i,a})b(g_{i,b})\exp(\lambda_b'g_{i,b})]$$

$$= \frac{E[a(g_{i,a})\exp(\lambda_a'g_{i,a})]}{E[\exp(\lambda_a'g_{i,a})]} \frac{E[b(g_{i,b})\exp(\lambda_b'g_{i,b})]}{E[\exp(\lambda_b'g_{i,b})]} = Q_a Q_b \equiv Q_{ab}$$

by the independence of $g_{i,a}$ and $g_{i,b}$ under the true untilted distribution. Hence $\text{plim}\,\hat{Q}_{ab} = \text{plim}\,\hat{Q}_a\,\text{plim}\,\hat{Q}_b$ as claimed. A similar result does not hold for EL because $(1-\lambda'g_i)^{-1} \neq (1-\lambda_a'g_{i,a})^{-1}(1-\lambda_b'g_{i,b})^{-1}$, unless $\max_{i\leq n}|\lambda'g_i| \xrightarrow{p} 0$, which is impossible under global misspecification. □

PROOF OF LEMMA 9. From (42), the first-order condition for $\hat{\theta}$ is

$$(52) \qquad -\tilde{G}'\tilde{\Omega}^{-1}\hat{g} - \sum_i \hat{w}_i G_i'\hat{\lambda}g_i'\tilde{\Omega}^{-1}\hat{g} + n^{-1}\sum_i G_i'\hat{\lambda} - \sum_i \hat{w}_i G_i'\hat{\lambda} = 0$$

(after transposition), where $\hat{\lambda}$ satisfies

$$(53) \qquad \sum_i \exp(\hat{\lambda}'g_i)g_i = 0.$$

Equation (52) contains products of sample moments which are difficult to analyze. Our goal is thus to define auxiliary parameters that will allow us to rewrite the first-order conditions as a linear function of sample moments.

Let us introduce the quantity $\hat{\tau}_i = \exp(\hat{\lambda}'g_i)$ and

$$(54) \qquad \hat{\tau} = n^{-1}\sum_i \hat{\tau}_i.$$

Noting that $\hat{w}_i = n^{-1}\hat{\tau}_i/\hat{\tau}$, (52) becomes

$$(55) \qquad \begin{aligned} &-\left(n^{-1}\sum_i \hat{\tau}_i G_i'\right)\left(n^{-1}\sum_i \hat{\tau}_i g_i g_i'\right)^{-1}\hat{g} \\ &- n^{-1}\sum_i \hat{\tau}_i G_i'\hat{\lambda}g_i'\left(n^{-1}\sum_i \hat{\tau}_i g_i g_i'\right)^{-1}\hat{g} \\ &+ n^{-1}\sum_i G_i'\hat{\lambda} - \frac{1}{\hat{\tau}}n^{-1}\sum_i \hat{\tau}_i G_i'\hat{\lambda} = 0. \end{aligned}$$

Now, we introduce $\hat{\kappa} = -(n^{-1}\sum_i(\hat{\tau}_i/\hat{\tau})g_i g_i')^{-1}\hat{g}$, or equivalently,

$$(56) \qquad \left(n^{-1}\sum_i \hat{\tau}_i g_i g_i'\right)\hat{\kappa} + \hat{\tau} n^{-1}\sum_i g_i = 0.$$

Substituting the $\hat{\kappa}$ whenever it appears in (55), after multiplying through by $\hat{\tau}$, yields

$$(57) \quad n^{-1}\sum_i \hat{\tau}_i G'_i \hat{\kappa} + n^{-1}\sum_i \hat{\tau}_i G'_i \hat{\lambda} g'_i \hat{\kappa} + n^{-1}\sum_i \hat{\tau} G'_i \hat{\lambda} - n^{-1}\sum_i \hat{\tau}_i G'_i \hat{\lambda} = 0.$$

Equation (57) is now linear in the sample moments. Equations (54), (53), (56) and (57) can be collected into a single vector of moment conditions $n^{-1}\sum_i \phi(x_i, \hat{\beta}) = 0$, where $\hat{\beta} = (\hat{\tau}, \hat{\kappa}', \hat{\lambda}', \hat{\theta}')'$ and

$$(58) \qquad \phi(x_i, \hat{\beta}) = \begin{bmatrix} \hat{\tau}_i - \hat{\tau} \\ \hat{\tau}_i g_i \\ (\hat{\tau} - \hat{\tau}_i) g_i + \hat{\tau}_i g_i g'_i \hat{\kappa} \\ \hat{\tau}_i G'_i \hat{\kappa} + \hat{\tau}_i G'_i \hat{\lambda} g'_i \hat{\kappa} - \hat{\tau}_i G'_i \hat{\lambda} + \hat{\tau} G'_i \hat{\lambda} \end{bmatrix}.$$

[For convenience, the third block is obtained by subtracting (53) from (56).] Noting that $\frac{\partial \hat{\tau}_i}{\partial \lambda} = \hat{\tau}_i g_i$, $\frac{\partial \hat{\tau}}{\partial \lambda} = 0$ and $\frac{\partial \hat{\tau}_i}{\partial \theta} = \hat{\tau}_i G'_i \lambda$, the first expression for $\phi(x_i, \hat{\beta})$ in (23) also follows. □

PROOF OF THEOREM 10. We first establish consistency of $\hat{\beta}$ in three steps: (i) Show that $\hat{\lambda}(\theta) \xrightarrow{p} \lambda^*(\theta)$ uniformly for $\theta \in \Theta$. (ii) Show that $\hat{\theta} \xrightarrow{p} \theta^*$ and therefore that $\hat{\lambda}(\hat{\theta}) \xrightarrow{p} \lambda^*(\theta^*)$. (iii) Show that this implies $\hat{\tau} \xrightarrow{p} \tau^*$ and $\hat{\kappa} \xrightarrow{p} \kappa^*$.

*Step* 1. By Lemma 2.4 in [46], continuity of $\exp(\lambda' g(x_i, \theta))$ in $\lambda$ and $\theta$, parts 1 and 4 of Assumption 3 imply that $\hat{M}_\theta(\lambda) \equiv n^{-1}\sum_i \exp(\lambda' g(x_i, \theta)) \xrightarrow{p} M_\theta(\lambda) \equiv E[\exp(\lambda' g(x_i, \theta))]$ uniformly over the compact set $\{(\lambda', \theta')' : \lambda \in \Lambda(\theta), \theta \in \Theta\}$, where $\Lambda(\theta)$ is as in part 4 of Assumption 3. We can then show that for any $\eta > 0$, $P[\sup_{\theta \in \Theta} \|\bar{\lambda}(\theta) - \lambda^*(\theta)\| \leq \eta] \to 1$, where $\bar{\lambda}(\theta) = \arg\min_{\lambda \in \Lambda(\theta)} \hat{M}_\theta(\lambda)$ as follows. For a given $\eta > 0$, select $\varepsilon = \inf_{\theta \in \Theta} \times \inf_{\lambda \in \Lambda(\theta): \|\lambda - \lambda^*(\theta)\| \geq \eta}(M_\theta(\lambda) - M_\theta(\lambda^*(\theta)))$, which is nonzero by the strict convexity of $M_\theta(\lambda)$ in $\lambda$ and the fact that $\Theta$ is compact. By the definition of $\varepsilon$, whenever $\sup_\theta (M_\theta(\bar{\lambda}(\theta)) - M_\theta(\lambda^*(\theta))) \leq \varepsilon$, then $\sup_{\theta \in \Theta} \|\bar{\lambda}(\theta) - \lambda^*(\theta)\| \leq \eta$. However, using the fact that $(\hat{M}_\theta(\bar{\lambda}(\theta)) - \hat{M}_\theta(\lambda^*(\theta))) < 0$, we have

$$\sup_\theta (M_\theta(\bar{\lambda}(\theta)) - M_\theta(\lambda^*(\theta)))$$
$$\leq \sup_\theta (M_\theta(\bar{\lambda}(\theta)) - \hat{M}_\theta(\bar{\lambda}(\theta))) + \sup_\theta (\hat{M}_\theta(\bar{\lambda}(\theta)) - \hat{M}_\theta(\lambda^*(\theta)))$$
$$\quad + \sup_\theta (\hat{M}_\theta(\lambda^*(\theta)) - M_\theta(\lambda^*(\theta)))$$
$$\leq \sup_\theta |M_\theta(\bar{\lambda}(\theta)) - \hat{M}_\theta(\bar{\lambda}(\theta))| + \sup_\theta |\hat{M}_\theta(\lambda^*(\theta)) - M_\theta(\lambda^*(\theta))|$$
$$\leq \frac{\varepsilon}{2} + \frac{\varepsilon}{2} = \varepsilon$$

w.p.a. 1. Hence, $\sup_{\theta\in\Theta}\|\bar{\lambda}(\theta)-\lambda^*(\theta)\|\leq \eta$ w.p.a. 1. In order to obtain the same conclusion for $\hat{\lambda}(\theta)$ rather than $\bar{\lambda}(\theta)$, we employ an argument similar to the proof of Theorem 2.7 in [46]. Since $\hat{M}_\theta(\lambda)$ is convex in $\lambda$ for any $\theta$, if the minimum $\bar{\lambda}(\theta)$ lies in the interior of $\Lambda(\theta)$, no other points in the complement of $\Lambda(\theta)$ can achieve a lower value and thus minimizing $\hat{M}_\theta(\lambda)$ over $\Lambda(\theta)$ or $\mathbb{R}^{N_g}$ yields the same answer asymptotically. This establishes $\sup_{\theta\in\Theta}\|\hat{\lambda}(\theta)-\lambda^*(\theta)\| \overset{p}{\to} 0$.

*Step* 2. $\ln \hat{L}(\theta) \equiv -\ln(n^{-1}\sum_i \exp(\hat{\lambda}'(\theta)g(x_i,\theta))) + \hat{\lambda}'(\theta)\hat{g}(\theta) \overset{p}{\to} \ln L(\theta)$ uniformly for $\theta\in\Theta$, because (i) $\sup_{\theta\in\Theta}\|\hat{\lambda}(\theta)-\lambda^*(\theta)\| \overset{p}{\to} 0$; (ii) $\sup_{\theta\in\Theta}\|\hat{g}(\theta) - E[g(x_i,\theta)]\| \overset{p}{\to} 0$ since $g(x_i,\theta)$ is continuous in $\theta$ and $E[\sup_{\theta\in\Theta}\|g(x_i,\theta)\|] < \infty$ by part 4 of Assumption 3 and by the inequality $|s| \leq \exp(-s)+\exp(s)$ for any $s\in\mathbb{R}$; and (iii) $\exp(\lambda'g(x_i,\theta))$ is continuous in $\theta$ and $E[\sup_{\theta\in\Theta}\sup_{\lambda\in\Lambda(\theta)} \exp(\lambda'g(x_i,\theta))] < \infty$ by part 4 of Assumption 3 (using Lemma 2.4 in [46]). Since $\ln L(\theta)$ is uniquely maximized at $\theta^*$, this implies, along with the uniform convergence of $\ln \hat{L}(\theta)$ and its continuity, that $\hat{\theta} \overset{p}{\to} \theta^*$. Since $\sup_{\theta\in\Theta}\|\hat{\lambda}(\theta)-\lambda^*(\theta)\| \overset{p}{\to} 0$ we also have that $\hat{\lambda}(\hat{\theta}) \overset{p}{\to} \lambda^*(\theta^*)$.

*Step* 3. As we have shown that $\hat{\theta} \overset{p}{\to} \theta^*$ and $\hat{\lambda} \overset{p}{\to} \lambda^*$ and since $\hat{\tau}$ and $\hat{\kappa}$ can be written as explicit continuous functions of $\hat{\lambda}$ and $\hat{\theta}$, by (54) and (56), it follows that $\hat{\tau} \overset{p}{\to} E[\tau_i] \equiv \tau^*$ and $\hat{\kappa} \overset{p}{\to} (E[\tau_i g_i(x_i,\theta^*)g_i'(x_i,\theta^*)])^{-1}(\tau^* E[g_i(x_i,\theta^*)]) \equiv \kappa^*$, where the fact that $E[\tau_i g_i(x_i,\theta^*)g_i'(x_i,\theta^*)]$ is invertible is implied by the assumption that $\Gamma$ is nonsingular.

Having established that $\hat{\beta} \overset{p}{\to} \beta^*$, we now turn to asymptotic normality. Since Lemma 9 defines a just-identified GMM estimator, we can use Theorem 3.4 in [46], specialized to the just-identified case, if we can show that (i) $E[\sup_{\beta\in\mathcal{B}}\|\partial\phi(x_i,\beta)/\partial\beta\|] < \infty$ for some neighborhood $\mathcal{B}$ of $\beta^*$ and that (ii) $E[\phi(x_i,\beta^*)\phi'(x_i,\beta^*)]$ exists.

The matrix $\partial\phi(x_i,\beta)/\partial\beta'$ consists of terms of the form $\alpha\exp(k_\tau\lambda'g_i)g^{k_g}G^{k_G}\times S^{k_S}$ for $0\leq k_g+k_G+k_S\leq 3$ and $k_\tau=0,1$, and where $g$, $G$, and $S$, respectively, denote elements of $g_i$, $G_i$, and $S_{jl}(x_i,\theta)$ and where $\alpha$ denotes products of elements of $\beta$ that are necessarily bounded for $\beta\in\mathcal{B}$. By part 6 of Assumption 3, we can establish (i): $\exp(k_\tau\lambda'g_i)|g|^{k_g}|G|^{k_G}|S|^{k_S} \leq \exp(k_\tau\lambda'g_i)\times |b(x_i)|^{k_g+k_G+k_S} \Rightarrow E[\sup_{\beta\in\mathcal{B}}\exp(k_\tau\lambda'g_i)|g|^{k_g}|G|^{k_G}|S|^{k_S}] \leq E[\sup_{\beta\in\mathcal{B}}\times \exp(k_\tau\lambda'g(x_i,\theta))(b(x_i))^{k_2}] = E[\sup_{\theta\in\mathcal{N}}\sup_{\lambda\in\Lambda(\theta)}\exp(k_\tau\lambda'g(x_i,\theta))(b(x_i))^{k_2}] < \infty$. The matrix $\phi(x_i,\beta)\phi_i'(x_i,\beta)$ has elements of the form $\alpha\exp(k_\tau\lambda'g_i)|g|^{k_g}|G|^{k_G}$ with $k_\tau=0,1,2$ and $0\leq k_g+k_G\leq 4$ and similar reasoning implies (ii). $\square$

40 S. M. SCHENNACH

Department of Economics
University of Chicago
1126 E. 59th Street
Chicago, Illinois 60637
USA
E-mail: smschenn@uchicago.edu